\documentclass{amsart}
\usepackage{amssymb}
\usepackage{amsthm}
\usepackage{thmtools}
\usepackage{mathrsfs}
\usepackage{stmaryrd}
\usepackage{enumitem}
\usepackage[french,english]{babel}
\usepackage[utf8]{inputenc}
\usepackage[T1]{fontenc}
\usepackage{lmodern}
\usepackage[all]{xy}
\usepackage[colorlinks=true]{hyperref}
\hypersetup{colorlinks, citecolor=blue, filecolor=black, linkcolor=black, urlcolor=black}
\usepackage{upref}
\usepackage{array}
\usepackage{color}
\usepackage{tikz-cd}
\usepackage{verbatim}
\usepackage{framed}
\usepackage[capitalize]{cleveref}

\usepackage[final]{showlabels}

\makeatletter\def\SL@eqnlefttext #1{\hbox to 0pt{\kern 60pt %or something else
\llap{\SL@margintext{#1}\quad}\hss}}
\makeatother

\def\mm{\mathfrak{m}}
\def\bC{\mathbf{C}}

\def\oo{\mathscr{O}_{\bC}}

\def\OK{\mathscr{O}_K}

\def\bQ{\mathbb{Q}}
\def\bQ0{\mathbb{Q}_{\ge 0}}

\def\Rep{\mathbf{Rep}}

\def\bfHig{\mathbf{Hig}}
\def\bfBun{\mathbf{Bun}}
\def\bfRep{\mathbf{Rep}}

\def\BB{\mathbb{B}}
\def\Gm{\mathbb{G}_m}

\def\bH{\mathbb{H}}

\def\Zp{\mathbb{Z}_p}
\def\Qp{\mathbb{Q}_p}
\def\BQ0{\mathbb{Q}_{\ge 0}}

\def\hatO{\widehat{\mathscr{O}}}

\def\dia{\lozenge}

\DeclareMathOperator{\GL}{GL}

\DeclareMathOperator{\ad}{ad}

\DeclareMathOperator{\Lie}{Lie}

\DeclareMathOperator{\HB}{HB}

\DeclareMathOperator{\gr}{gr}

\DeclareMathOperator{\can}{can}

\DeclareMathOperator{\et}{\textnormal{\'{e}t}}
\DeclareMathOperator{\proet}{\textnormal{pro\'{e}t}}

\DeclareMathOperator{\Gal}{Gal}
\DeclareMathOperator{\Ker}{Ker}

\DeclareMathOperator{\Hom}{Hom}

\DeclareMathOperator{\Sym}{Sym}

\DeclareMathOperator{\Ext}{Ext}

\DeclareMathOperator{\Spec}{Spec}

\DeclareMathOperator{\Spd}{Spd}
\DeclareMathOperator{\Spa}{Spa}

\DeclareMathOperator{\cont}{cont}
\DeclareMathOperator{\rH}{H}

\DeclareMathOperator{\rB}{B}

\DeclareMathOperator{\rS}{S}

\DeclareMathOperator{\Gr}{Gr}

\DeclareMathOperator{\id}{id}

\DeclareMathOperator{\dR}{dR}

\DeclareMathOperator{\Exp}{Exp}

\DeclareMathOperator{\Bun}{Bun}

\DeclareMathOperator{\HT}{HT}

\DeclareMathOperator{\Sec}{Sec}

\DeclareMathOperator{\Perf}{Perf}
\DeclareMathOperator{\univ}{univ}

\newtheorem{theorem}{Theorem}[subsection]
\newtheorem{prop}[theorem]{Proposition}
\newtheorem{lemma}[theorem]{Lemma}

\newtheorem{coro}[theorem]{Corollary}

\theoremstyle{definition}
\newtheorem{rem}[theorem]{Remark}

\newtheorem{secnumber}[theorem]{}
\newtheorem{question}[theorem]{Question}

\crefname{theorem}{theorem}{theorems}
\Crefname{theorem}{Theorem}{Theorems}

\crefname{prop}{proposition}{propositions}
\Crefname{prop}{Proposition}{Propositions}

\crefname{lemma}{lemma}{lemmas}
\Crefname{lemma}{Lemma}{Lemmas}

\crefname{coro}{corollary}{corollaries}
\Crefname{coro}{Corollary}{Corollaries}

\crefname{rem}{remark}{remarks}
\Crefname{rem}{Remark}{Remarks}

\crefname{definition}{definition}{definitions}
\Crefname{definition}{Definition}{Definitions}

\crefname{conjecture}{conjecture}{conjectures}
\Crefname{conjecture}{Conjecture}{Conjectures}

\crefname{question}{question}{questions}
\Crefname{question}{Question}{Questions}

\crefname{secnumber}{paragraph}{paragraphs}
\Crefname{secnumber}{Paragraph}{Paragraphs}

\numberwithin{equation}{section}
\numberwithin{equation}{theorem}

\makeatletter
\newsavebox{\@brx}
\newcommand{\llangle}[1][]{\savebox{\@brx}{\(\m@th{#1\langle}\)}%
  \mathopen{\copy\@brx\kern-0.5\wd\@brx\usebox{\@brx}}}
\newcommand{\rrangle}[1][]{\savebox{\@brx}{\(\m@th{#1\rangle}\)}%
  \mathclose{\copy\@brx\kern-0.5\wd\@brx\usebox{\@brx}}}
\makeatother

\newcommand{\quash}[1]{}

\AtEndDocument{\bigskip{\footnotesize

\textsc{Daxin Xu, Morningside Center of Mathematics and Hua Loo-Keng Key Laboratory of Mathematics, Academy of Mathematics and Systems Science, Chinese Academy of Sciences, Beijing 100190, China} \par
\textit{E-mail address}: \texttt{daxin.xu@amss.ac.cn} \par
  \textsc{Xinwen Zhu, Department of Mathematics, Building 380, Stanford, California 94305} \par
  \textit{E-mail address}: \texttt{zhuxw@stanford.edu} \par

}}

\title{$p$-adic representations for rank-two nilpotent Higgs bundles}
\author{Daxin Xu, Xinwen Zhu}
\date{\today}

\begin{document}
\selectlanguage{english}

\begin{abstract}
	Let $C$ be a smooth projective $p$-adic curve of genus $\ge 2$.
	We show that every rank-two semi-stable Higgs bundle of degree zero on $C$ with non-zero nilpotent Higgs field corresponds to a continuous $p$-adic representation of the geometric fundamental group of $C$ under the $p$-adic Simpson correspondence associated with a suitable lifting of $C$. 
	For the representation $\rho$ associated with the uniformizing Higgs bundle, we construct a Hodge--Tate map from the pro-finite-\'etale cover of $C$, trivializing $\rho$, to $\mathbb{P}^1$. 
\end{abstract}

\maketitle
\section{Introduction} \label{s:intro}

\subsection{Background on the complex and $p$-adic Simpson correspondence}

\begin{secnumber}
	Let $X$ be a compact Riemann surface of genus $g\ge 2$.
	The classical Narasimhan--Seshadri theorem \cite{NS65} establishes a bijection between stable holomorphic vector bundles of degree zero on $X$ and irreducible unitary representations of the topological fundamental group $\pi_1(X)$. 
	It was generalized to higher-dimensional compact K\"ahler manifolds by Donaldson \cite{Don85} and Uhlenbeck--Yau \cite{UY86} via the existence of Hermitian--Einstein metrics. 
	In his seminal work \cite{Hit87}, Hitchin extended this picture to non-unitary representations by introducing Higgs bundles $(E,\theta)$ on $X$ and studying the self-duality equations for a Hermitian metric on $E$.
	Subsequently, Simpson \cite{Simpson} generalized Hitchin's theorem to higher-dimensional smooth projective complex varieties, establishing the non-abelian Hodge correspondence between irreducible finite-dimensional complex representations of $\pi_1^{\textnormal{top}}(X)$ and stable Higgs bundles with vanishing Chern classes.

	A fundamental example in Hitchin's work \cite[\S~11]{Hit87} is the \textit{uniformizing Higgs bundle}. 
	Choosing a theta characteristic $L=\Omega_X^{1/2}$ on $X$, i.e.\ a holomorphic line bundle equipped with an isomorphism $L^{\otimes 2}\xrightarrow{\sim}\Omega_X^1$, one considers the rank-two nilpotent Higgs bundle
	\begin{equation} \label{eq:uniformizing-Higgs}
	\left(E=L^{-1}\oplus L,\; \Theta=\begin{pmatrix} 0 & \id \\ 0 & 0 \end{pmatrix}: E\longrightarrow E\otimes_{\mathscr{O}_X}\Omega_X^1\right).
	\end{equation}
	Although the underlying vector bundle $L^{-1}\oplus L$ is unstable, the Higgs bundle $(E,\Theta)$ is stable.
	Hitchin showed that the associated projective representation $\rho:\pi_1^{\textnormal{top}}(X)\to \mathrm{PGL}_2(\mathbb{C})$ factors through $\mathrm{PSL}_2(\mathbb{R})$ and realize the complex uniformization of $X$ by the upper half-plane $\mathfrak{H}\subset \mathbb{P}^1(\mathbb{C})$. 
	\[
	\xymatrix{ 
	\mathfrak{H} \ar[r] \ar[d]
	   & \mathbb{P}^1(\mathbb{C}) & \tau\mapsto [\tau:1] \\
	 X & &
}
	\]

	The main goal of the present article is to study a $p$-adic analogue of Hitchin's uniformization result. 
	Namely, we construct a $p$-adic representation $\rho$ of the geometric fundamental group associated with the uniformizing Higgs bundle %(and, more generally, rank-two semi-stable Higgs bundles of degree zero with non-zero nilpotent Higgs field) 
	under the $p$-adic Simpson correspondence over a $p$-adic curve $C$ of genus $\ge 2$ defined by a suitable lifting of $C$. 
	This result also applies to 
	rank-two semi-stable Higgs bundles of degree zero with non-zero nilpotent Higgs field. 
	The minimal pro-finite-\'etale cover $\widetilde{C}_{\rho}$ over $C$, trivializing $\rho$, is pointwise perfectoid in the sense of \cite{HeSen}. We construct a proper, surjective and equivariant Hodge--Tate map $\pi_{\HT,\rho}:\widetilde{C}_{\rho}\to \mathbb{P}^{1,\dia}$ of diamonds classifying the Hodge--Tate filtration on the pro-\'etale vector bundle associated with $\rho$. 
	\[
	\xymatrix{
	\widetilde{C}_{\rho} \ar[r]^{\pi_{\HT,\rho}} \ar[d]_{\pi}
	   & \mathbb{P}^{1,\dia} \\
	 C &
	}
	\]
	This construction provides a $p$-adic period-map analogue of complex uniformization. 
	It is compatible with Scholze's Hodge--Tate period map \cite{Sch-HT} for Shimura curves (and for modular curves in the logarithmic setting).
\end{secnumber}

\begin{secnumber}
	Over $p$-adic smooth proper varieties, Faltings \cite{Fal05} initiated an analogue of Simpson's theory, known as the \textit{$p$-adic Simpson correspondence}.
	Let $K$ be a complete discrete valuation field of characteristic zero with residue field $k=\overline{\mathbb{F}}_p$ of characteristic $p>0$, $\OK$ its valuation ring, $\overline{K}$ an algebraic closure of $K$, and $\bC$ (resp.\ $\oo$) the $p$-adic completion of $\overline{K}$ (resp.\ $\mathscr{O}_{\overline{K}}$).
	Let $X$ be a connected smooth proper rigid analytic variety over $\bC$.
	In modern terms, the global $p$-adic Simpson correspondence establishes an equivalence of categories between pro-\'etale vector bundles on $X_{\proet}$ and Higgs bundles on $X$:
	\begin{equation}\label{eq:intro-v-Simpson}
		\rH_{\mathcal{X},\Exp}: \Bun_{\proet}(X)\xrightarrow{\;\sim\;}\HB(X).
	\end{equation}
	This equivalence depends on the choices of a flat lifting $\mathcal{X}$ of $X$ over $\rB_{\dR,2}^+=\rB_{\dR}^+/\xi^2\rB_{\dR}^+$ and an exponential map $\Exp:(\bC,+)\to (1+\mm_{\bC},\times)$ splitting the $p$-adic logarithm $\log:1+\mm_{\bC}\to \bC$.

	There are mainly two approaches to the equivalence \eqref{eq:intro-v-Simpson}:
	\begin{enumerate}
		\item[\textnormal{(i)}] The original approach is introduced by Faltings \cite{Fal05} and developed by Abbes--Gros, Tsuji \cite{AGT, AbbesGros, AGrel}, and Wang \cite{Wang23}. It consists of two steps:
		(a) establishing a ``small'' equivalence between objects on both sides that are $p$-adically close to the trivial object;
		(b) obtaining the full equivalence by twisted pullback and \'etale descent (this step was only estabilished for smooth proper curves).

	\item[\textnormal{(ii)}] Another approach is based on the Hitchin fibration and the geometric Sen operator (defined by the ``small'' equivalence), recently established by Heuer \cite{Heu25} for smooth proper rigid analytic varieties over $\bC$, and geometrized in the curve case via moduli stacks by Heuer--Xu \cite{HX24} (see also the related work of Bhatt--Zhang via the Simpson gerbe \cite{Zha25}). 
	\end{enumerate}
\end{secnumber}

\begin{secnumber}
	Fix a geometric point $\overline{x}$ of $X$.
	There is a natural fully faithful functor sending a continuous representation of the \'etale fundamental group to its associated pro-\'etale vector bundle:
	\[
	\Rep(\pi_1(X,\overline{x}),\bC)\hookrightarrow \Bun_{\proet}(X),\qquad
	\rho\longmapsto V_\rho.
	\]
	Composing this embedding with \eqref{eq:intro-v-Simpson} yields a fully faithful functor
	\[
	\bH_{\mathcal{X},\Exp}: \Rep(\pi_1(X,\overline{x}),\bC)\hookrightarrow \HB(X).
	\]
	A fundamental problem in the $p$-adic Simpson correspondence is to characterize which Higgs bundles lie in the essential image of $\bH_{\mathcal{X},\Exp}$. 
	
	When $C$ is a smooth projective connected curve over $\bC$, every Higgs bundle in the essential image of $\bH_{\mathcal{C},\Exp}$ is semi-stable of degree zero. Faltings first sketched a proof of this statement in \cite[\S~5]{Fal05}; see also \cite{Xu22,HX24}. 
	Inspired by the complex Simpson correspondence, one may formulate the following question:
\end{secnumber}

\begin{question}\label{q:intro}
%	Which Higgs bundles over $C$ lie in the essential image of $\bH_{\mathcal{C},\Exp}$?
	Does every semi-stable Higgs bundle of degree zero over $C$ lie in the essential image of $\bH_{\mathcal{C},\Exp}$?
\end{question}

	A positive answer to the above question is known for Higgs line bundles of degree zero \cite{DW05, Heu24, Xu22}.
	In higher rank, the question remains widely open, even in the case of vector bundles with zero Higgs field. 

	In \cite{Xu22}, we characterized this essential image for curves of genus at least two defined over $\overline K$, based on Deninger--Werner's parallel transport functor \cite{DW05}. 
	The latter is compatible with the Higgs-field zero case of the $p$-adic Simpson correspondence \cite{Xu17}.

\subsection{Main results: rank two nilpotent Higgs bundles and the Hodge--Tate map}

	Our first result answers Question~\ref{q:intro} for nilpotent Higgs bundles of rank two with non-zero Higgs field, allowing the lifting to vary:

\begin{theorem}\label{t:intro-nilpotent}
	Let $C$ be a smooth projective curve of genus $g\ge 2$, and let $(E,\Theta)$ be a rank-two semi-stable Higgs bundle of degree zero with non-zero nilpotent Higgs field $\Theta$ over $C$.
	There exists a flat lifting $\mathcal{C}$ of $C$ over $\rB_{\dR,2}^+$, such that $(E,\Theta)$ lies in the essential image of
	\[
	\bH_{\mathcal{C},\Exp}: \Rep(\pi_1(C,\overline{x}),\bC)\longrightarrow \HB(C)
	\]
	for any choice of $\Exp$ (and the associated representation of $\pi_1(C,\overline{x})$ is independent of $\Exp$). 
\end{theorem}

We outline the proof of \Cref{t:intro-nilpotent} (\S~\ref{s:main-result}).
	After twisting by a line bundle of degree zero, we may assume $\det(E,\Theta)\simeq (\mathscr{O}_C,0)$.
	Set $L^{-1}=\Ker(\Theta)$. Then $\deg L\ge 0$, and $(E,\Theta)$ is an extension of $(L,0)$ by $(L^{-1},0)$ determined by $\gamma\in \rH^1(C,L^{-2})$ and $0\neq \eta\in \rH^0(C,L^{-2}\otimes_{\mathscr{O}_C}\Omega_C^1)$, denoted by $E(\gamma,\eta)$.
	In \S~\ref{s:extensions}, we find a class $\gamma_0\in \rH^1(C,L^{-2})$ whose underlying vector bundle is an extension of degree-zero line bundles.
	For a fixed lifting $\mathcal{C}_0$, the Higgs bundle $E(\gamma_0,0)$ lies in the essential image of $\bH_{\mathcal{C}_0,\Exp}$ by \cite{Heu24,Heu25}. Then by the openness of the representation locus \cite{HX24}, so does $E(\gamma_0,a\eta)$ for a sufficiently small $a\in\bC^{\times}$. 
	By compatibility of $\rH_{\mathcal{C},\Exp}$ with twisted pullback (\S~\ref{s:complement}), changing $\mathcal{C}_0$ to another lifting $\mathcal{C}$ allows us to shift $\gamma_0$ to $a\gamma$. Then we deduce that $(E,\Theta)\simeq E(a\gamma,a\eta)$ lies in the essential image of $\bH_{\mathcal{C},\Exp}$. 

	In particular, \Cref{t:intro-nilpotent} applies to Hitchin's uniformizing Higgs bundle \eqref{eq:uniformizing-Higgs}. 
	By taking symmetric powers and the fact that $\bH_{\mathcal{C},\Exp}$ is a tensor functor, we obtain the following higher-rank consequence:

\begin{coro}\label{c:kostant}
	Let $L$ be a theta characteristic on $C$ equipped with $\theta:L^{\otimes 2}\xrightarrow{\sim}\Omega_C^1$.
	For $n\ge 2$, set $E_n=\bigoplus_{i=0}^{n-1}L^{1-n+2i}$ on $C$, and let $\Theta_n:E_n\to E_n\otimes_{\mathscr{O}_C}\Omega_C^1$ vanish on $L^{1-n}$ and send each remaining summand $L^a$ to $L^{a-2}\otimes_{\mathscr{O}_C}\Omega_C^1$ via the isomorphism induced by $\theta$.
	There exists a flat lifting $\mathcal{C}$ of $C$ over $\rB_{\dR,2}^+$ such that $(E_n,\Theta_n)$ lies in the essential image of $\bH_{\mathcal{C},\Exp}$ for every $n\ge 2$ and every choice of $\Exp$.
\end{coro}

\begin{secnumber}
	In \S~\ref{s:HT-map}, we apply Theorem~\ref{t:intro-nilpotent} to study a $p$-adic analogue of complex uniformization.
	Let $L$ be a theta characteristic on $C$ equipped with an isomorphism $L^{\otimes 2}\xrightarrow{\sim}\Omega_C^1$. 
	The uniformizing Higgs bundle $(E=L^{-1}\oplus L,\Theta)$ admits a filtration:
	\begin{equation} \label{eq:filtration-Higgs}
	0\to (L^{-1},0)\to (E,\Theta)\to (L,0)\to 0.
	\end{equation}
	
	By Theorem~\ref{t:intro-nilpotent}, there exists a flat lifting $\mathcal{C}$ of $C$ over $\rB_{\dR,2}^+$ such that $(E,\Theta)$ corresponds under $\rH_{\mathcal{C},\Exp}$ to a continuous representation
	\[
	\rho:\pi_1(C,\overline{x})\longrightarrow \GL_2(\bC).
	\]
	Fixing a $\rho$-stable $\oo$-lattice in $\bC^{\oplus 2}$, let $C_n\to C$ be the connected finite \'etale cover trivializing $\rho\bmod p^n$, and consider the minimal pro-finite-\'etale cover in the category of diamonds over $\Spd(\bC,\oo)$:
	\[
	\pi:\widetilde{C}_{\rho}=\varprojlim_{n\ge 1} C_n \longrightarrow C,
	\]
	which trivializes $\rho$ and its associated pro-\'etale vector bundle $V_\rho$.
	Under the $p$-adic Simpson correspondence, we obtain a short exact sequence of pro-\'etale vector bundles from \eqref{eq:filtration-Higgs}:
	\[0\to \nu^*L^{-1}\to V_\rho\to \nu^*L\to 0,\]
	where $\nu:C_{\proet}\to C_{\et}$ is the canonical morphism. 
\end{secnumber}

\begin{theorem}[Theorem~\ref{t:HT-map}]\label{t:intro-HT}
	\textnormal{(i)} Assume in addition that $C$ descends to a smooth projective curve over
    $\overline K$. The pro-finite-\'etale cover $\widetilde{C}_{\rho}$ is pointwise perfectoid in the sense of \cite{HeSen}.

    \textnormal{(ii)} There exists a proper and surjective $\pi_1(C,\overline{x})$-equivariant Hodge--Tate map of diamonds over $\Spd(\bC,\oo)$:
	\[
	\pi_{\HT,\rho}:\widetilde{C}_{\rho}\longrightarrow \mathbb{P}^{1,\dia}_{\bC}
	\]
	classifying the rank-one quotient $\pi^*V_\rho\simeq \widehat{\mathscr{O}}^{\oplus 2}\twoheadrightarrow \pi^*\nu^*L$.
	In particular, we have 
	\begin{equation} \label{eq:pullback-line-bundles}
		\pi_{\HT,\rho}^*(\mathscr{O}_{\mathbb{P}^1}(1))\simeq \pi^*L.
\end{equation}
\end{theorem}
%\dx{we may define it from the projective representation of $\rho$ so that $\pi_{\HT}$ depends only on the lift $\mathcal{C}$.}

\Cref{t:intro-nilpotent} is know for Shimura (resp. modular) curves with the canonical lift via the $p$-adic Riemann--Hilbert correspondence \cite{LZ17,DLLZ23}, see \S~\ref{ss:comments}. 
	In this case, \Cref{t:intro-HT} (resp. its logarithmic variant) is compatible with Scholze's Hodge--Tate period map \cite{Sch-HT} and the cover $\widetilde{C}_{\rho}$ is known to be perfectoid. 

	The pro-\'etale bundle $V_\rho$ is isomorphic the Faltings extension defined by $\mathcal{C}$ up to twisting $L$ and is trivialized by $\pi$. 
	Then the pointwise perfectoidness (i) follows from He's Sen-theoretic criterion \cite[Theorem 7.19]{HeSen}. 

	For Theorem~\ref{t:intro-HT}(ii), the properness follows from that of $\widetilde{C}_{\rho}\to \Spd(\bC,\oo)$. 
	Therefore, it suffices to show the surjectivity on classical points. 
	If $\pi_{\HT,\rho}$ factors through $\mathbb{A}^{1,\dia}_{\bC}$, then $\pi_{\HT,\rho}^*(\mathscr{O}_{\mathbb{P}^1}(1))$ should be a trivial bundle on $\widetilde{C}_{\rho}$. 
	This would contradict to \eqref{eq:pullback-line-bundles} and that $\deg L=g-1>0$. 

\subsection{Structure of the article}

	In \S~\ref{s:complement}, we establish the compatibility of the $p$-adic Simpson correspondence $\rH_{\mathcal{C},\Exp}$ for proper curves \cite{Heu25,HX24} with the twisted pullback functor \cite{Fal05,Xu22}.
	We also compare the $\rH^1(C,\xi T_C)$-torsors of flat liftings, Faltings extensions and rank-two nilpotent Higgs extensions.
	In \S~\ref{s:extensions}, we use higher secant varieties of curves to construct semi-stable, non-stable extensions of line bundles, which are extensions of line bundles of degree zero. 
	Section~\ref{s:main-result} is devoted to the proof of Theorem~\ref{t:intro-nilpotent}, together with a discussion on an alternative proof via \cite{Xu22} and on the canonical lifting. 
	Finally, in \S~\ref{s:HT-map}, we prove Theorem~\ref{t:intro-HT} on the pointwise perfectoidness of $\widetilde{C}_{\rho}$ and the Hodge--Tate map $\pi_{\HT,\rho}$. \\

\textbf{Acknowledgements.} 
	We would like to thank Ahmed Abbes, Tongmu He, Ben Heuer, Sasha Petrov, Yupeng Wang and Kang Zuo for valuable discussions and comments. 
	D. Xu is supported by National Key R\&D Program of China
2025YFA1018000, National Natural Science Foundation of China Grant no. 12625103, 12288201 and CAS Project for Young Scientists in Basic Research, Grant No. YSBR-033.
	 X. Zhu is supported by NSF grant DMS-2200940.

\textbf{AI disclosure.}
This project began in 2025 and was completed largely without the use of AI. We used large language models to assist with exposition, literature search, and proofreading. All mathematical results and proofs are our own, and we take full responsibility for the content of this article.

\section{Complement on the $p$-adic Simpson correspondence}
\label{s:complement}	

\subsection{The $p$-adic Simpson correspondence over proper curves and twisted pullback}
\label{ss:comparison}

Let $C$ be a smooth projective $\bC$-curve of genus $\geq2$. 
We use the same notation for its rigid analytification.  
We denote the canonical morphism of ringed topoi by 
\[\nu:(C_{\proet},\widehat{\mathscr{O}}_C)\to (C_{\et},\mathscr{O}_C).\]
Fix a flat lifting $\mathcal C$ of $C$ over $\rB_{\dR,2}^+$ and an
exponential $\Exp:\bC\to1+\mathfrak m_{\bC}$.
In this subsection, we show the compatibility of the $p$-adic Simpson correspondence over proper curves \cite{Heu25,HX24} with the twisted pullback functor in \cite{Fal05,Xu22}.

\begin{secnumber} \label{sss:review-LExp}
We first recall the twisted pullback functor $f^{\circ}_{\mathcal{C}',\mathcal{C},\Exp}$ and the construction of the twisting \'etale line bundle $\mathcal L_{f,b}^{\Exp}$ via rigidification from \cite[Appendix~A.2]{Xu22}.
Let $f:C'\to C$ be a finite morphism of smooth projective $\bC$-curves of genus $\ge 2$, and let $\mathcal C', \mathcal{C}$ be flat liftings of $C',C$ over $\rB_{\dR,2}^+$, respectively.
For a point of the Hitchin base
\[
b=(b_i)_{i=1}^n \in A=\bigoplus_{i=1}^n\Gamma(C,\xi^{-i}\Omega_C^{\otimes i}),
\]
the associated spectral cover over $C'$ is
\begin{equation}\label{eq:C-fb}
 \pi_{f,b}:C'_{f,b}=\Spec_{\mathscr{O}_{C'}}(\mathcal{A}_{f,b})\longrightarrow C',
 \qquad \mathcal{A}_{f,b}=\Sym^{\bullet}_{\mathscr{O}_{C'}}(\xi f^*T_C)/\mathcal{I}_{f,b};
\end{equation}
see \cite[\S~A.2.1]{Xu22}.
Let $\mathcal L_f$ be the $\xi f^*T_C$-torsor on $C'_{\et}$ of local liftings
$\widetilde f:\mathcal C'\to\mathcal C$ of $f:C'\to C$. 

Choose $x\in C(\bC)$ with $b_n(x)=0$ and $x'\in C'(\bC)$ above $x$.
Such a point exists since $b_n$ is a section of a positive-degree line bundle $\xi^{-n}\Omega_C^{\otimes n}$. 
Then $x'$ has a unique lift $z_f=(x',0)\in C'_{f,b}$ such that if we set
\[
I:=\Ker(\mathcal{A}_{f,b}\to i_{z_f,*}\bC),
\]
the canonical map $\beta: \xi f^*T_C\to \mathcal{A}_{f,b}$ factors through $I$ \cite[Lemma~A.2.2]{Xu22}. 
Then $\beta_* \mathcal{L}_f$ defines an $I$-torsor on $C'$ up to a unique isomorphism. 

Let $\mathbf P_{\et,f,b}^{z_f}$ be the analytification of the algebraic
Picard group of $C'_{f,b}$ rigidified at $z_f$ \cite[(A.2.3.1)]{Xu22}.
Its Lie algebra is $\rH^1(C',I)$, and the exponential map
\cite[Proposition~A.1.2]{Xu22} defines an object in the neutral component of $\mathbf P_{\et,f,b}^{z_f}$:
\[
[ (\mathcal L_{f,b}^{\Exp},\alpha_{z_f}) ]
 =\Exp_{\mathbf P_{\et,f,b}^{z_f}}
 \bigl([\beta_* \mathcal{L}_f]\bigr) \in (\mathbf P_{\et,f,b}^{z_f})^{\circ}(\bC). 
\]
Since
$\mathcal A_{f,b}\simeq\bigoplus_{j=0}^{n-1}(\xi f^*T_C)^{\otimes j}$
and $\deg(f^*T_C)<0$, we have $\rH^0(C',\mathcal A_{f,b})=\bC$. 
An automorphism of $\mathcal L_{f,b}^{\Exp}$ preserving the rigidification is therefore the identity.
The line bundle $\mathcal L_{f,b}^{\Exp}$ is independent of $z_f$, up to a unique isomorphism
\cite[Proposition~A.2.4(i)]{Xu22}.

Via $f^*\Omega_C\to\Omega_{C'}$, tensoring with $\mathcal{L}_{f,b}^{\Exp}$ defines the twisted pullback functor \cite[Proposition~A.2.5]{Xu22}:
\[
 f^{\circ}_{\mathcal{C}',\mathcal{C},\Exp}:\HB(C)\longrightarrow \HB(C'),\qquad
 (E,\theta)\longmapsto f^*(E,\theta)\otimes_{\mathcal{A}_{f,b}}\mathcal{L}_{f,b}^{\Exp}.
\]
On nilpotent Higgs bundles ($b=0$), the augmentation ideal of $\mathcal{A}_{f,0}$ is nilpotent, so $\mathcal{L}_{f,0}^{\Exp}$ is obtained via the truncated exponential series and is independent of $\Exp$.
\end{secnumber}
\begin{secnumber} \label{sss:Higgs-Tate-torsor}
Let $\mathcal{C}$ be a flat lift of $C$ over $\rB_{\dR,2}^+$. 
We recall the Higgs--Tate torsor $\mathcal{P}_{\mathcal{C}}$ associated with $\mathcal{C}$ \cite[II.10.3]{AGT}, \cite[Definition 2.12]{Heu25}.
Let $U=\Spa(R,R^+)$ be an affinoid perfectoid object of $C_{\proet}$, and set $\BB_{\dR,2}^+(R)=\BB_{\dR}^+(R)/\xi^2$.
Choose an affinoid open subset $V=\Spa(A,A^+)\subset C$ through which $U\to C$ factors, and put $\widetilde{A}=\Gamma(V,\mathscr{O}_{\mathcal{C}})$, a flat $\rB_{\dR,2}^+$-algebra lifting $A$. 
Let $\mathcal{P}_{\mathcal{C}}(U)$ be the set of continuous $\rB_{\dR,2}^+$-algebra homomorphisms $\varphi:\widetilde{A}\to \rB_{\dR,2}^+(R)$ such that $\theta\circ\varphi$ is the composition $\widetilde{A}\to A\to R$, i.e. the set of liftings of $A\to R$ to $\rB_{\dR,2}^+$-algebra homomorphisms $\widetilde{A}\to \rB_{\dR,2}^+(R)$. 
Then $U\mapsto \mathcal{P}_{\mathcal{C}}(U)$ defines a torsor $\mathcal{P}_{\mathcal{C}}$ under $\nu^*(\xi T_C)$ on $C_{\proet}$ \cite[Lemma 2.13]{Heu25}. 

For a vector bundle $W$ on $C$, we have $R^1\nu_*\nu^*W\simeq W\otimes_{\mathscr{O}_C}\xi^{-1}\Omega_C$, by the projection formula and \cite[Proposition~3.23]{Sch13}.
Since $\rH^2(C,W)=0$, the Leray spectral sequence for $\nu$ gives an exact sequence 
\begin{equation}\label{eq:HT-sequence-W}
 0\longrightarrow\rH^1(C,W)\xrightarrow{\nu^*}
 \rH^1(C_{\proet},\nu^*W)\xrightarrow{\HT}
 \rH^0(C,W\otimes_{\mathscr{O}_C}\xi^{-1}\Omega_C)\longrightarrow0.
\end{equation}
For $W=\xi T_C$, the class $[\mathcal{P}_{\mathcal{C}}]\in \rH^1(C_{\proet},\nu^*(\xi T_C))$ maps to $1\in\rH^0(C,\mathscr{O}_C)$ under the \textit{Hodge--Tate map} $\HT$ \cite[Proposition 2.15]{Heu25}.
\end{secnumber}

\begin{secnumber}\label{sss:spectral-twist}
	Let $\mathcal{A}_b=\mathcal{A}_{\id_C,b}$ be the spectral algebra of \ref{sss:review-LExp}, $\pi_b:C_b=\Spec(\mathcal{A}_b)\to C$ the spectral cover.
Let $\beta_b:\xi T_C\to \mathcal{A}_b$ be the canonical map and $\tau_b\in \Gamma(C,\mathcal{A}_b\otimes \xi^{-1}\Omega_C)$ the associated tautological section. 

For $\mathcal A_b$, the multiplicative Hodge--Tate
sequence \cite[Proposition 2.9, Corollary 2.11 and Theorem 2.4]{Heu25}
defines the \textit{Hodge--Tate logarithm}:
\[
 \mathrm{HTlog}:\rH^1(C_{\proet},(\nu^*\mathcal A_b)^\times)
 \longrightarrow\Gamma(C,\mathcal A_b\otimes\xi^{-1}\Omega_C).
\]
whose kernel consists of the classes of \'etale invertible $\mathcal A_b$-modules. 

We define $\mathscr{H}_b$ as the groupoid of invertible $\nu^*\mathcal{A}_b$-modules $\mathcal{K}$ on $C_{\proet}$ such that 
\[
\mathrm{HTlog}(\mathcal{K})=\tau_b \in \Gamma(C,\mathcal{A}_b\otimes_{\mathscr{O}_C}\xi^{-1}\Omega_C),
\]
with $\nu^*\mathcal{A}_b$-linear isomorphisms as morphisms.
This is a slight variant of the twisting gerbe of \cite[\S~6.1]{HX24}. 
The Picard groupoid $\mathscr{P}_b$ of \'etale line bundles on $C_b$ acts on $\mathscr{H}_b$ by $(\mathcal{L},\mathcal{K})\mapsto \mathcal{K}\otimes_{\nu^*\mathcal{A}_b}\nu^*\mathcal{L}$. 
\end{secnumber}

\begin{secnumber}\label{sss:twisting-object}
Pushing out $\mathcal{P}_{\mathcal{C}}$ along $\nu^*\beta_b$, we obtain a $\nu^*\mathcal{A}_b$-torsor on $C_{\proet}$
\begin{equation}\label{eq:spectral-HT-torsor}
 \mathcal{P}_{\mathcal{C},b}=\mathcal{P}_{\mathcal{C}}\times^{\nu^*(\xi T_C),\,\nu^*\beta_b}\nu^*\mathcal{A}_b.
\end{equation}
By \ref{sss:Higgs-Tate-torsor}, the Hodge--Tate image of $[\mathcal{P}_{\mathcal{C},b}]$ in $\rH^0(C,\mathcal{A}_b\otimes_{\mathscr{O}_C}\xi^{-1}\Omega_C)$ is $(\beta_b\otimes\id)(1)=\tau_b$.

Next, we canonically construct an object $\mathcal{K}^{\Exp}_{\mathcal{C},b}$ of $\mathscr{H}_b$ from $\mathcal{P}_{\mathcal{C}}$, following the rigidified Picard construction of \cite[\S~3.2]{Heu25}, with the base point used in \ref{sss:review-LExp} \cite[\S~A.2.3]{Xu22}. 
We take a point $x\in C(\bC)$ such that $b_n(x)=0$. 
This defines a closed point of  $ z=(x,0)\in C_b$,
and an evaluation homomorphism 
\[
 \varepsilon_z:\mathcal A_b|_x\longrightarrow\bC
\]
such that $\varepsilon_z\circ\beta_b|_x=0$.
Then the restriction of $\mathcal{P}_{\mathcal{C},b}$ along $\varepsilon_z$ admits a
canonical trivialization $\alpha_{\can}$: the pushout of
$\mathcal P_{\mathcal C}|_x$ along the zero map
$\xi T_C|_x\to\bC$ is canonically the trivial additive $\bC$-torsor. 

Let $\mathbf P_{\proet,b}$ be the relative pro-\'etale Picard functor of
$\nu^*\mathcal A_b$ \cite[Definition~2.2]{Heu25}. It is represented by a rigid group \cite[Theorem~2.4(3)]{Heu25}. Let $\mathbf P_{\proet,b}^{z}$ be its rigidified version,
parametrizing invertible $\nu^*\mathcal A_b$-modules $\mathcal L$ equipped with a rigidification
\[
 \alpha_z:
 \mathcal L|_x\otimes_{\mathcal A_b|_x,\varepsilon_z}\bC
 \xrightarrow{\sim}\bC.
\]

Here is a comparison between $\mathbf P_{\proet,b}^{z}$ and the rigidified construction in \cite[Definition 3.12]{Heu25}.
Since
\[
 \mathcal A_b\simeq\bigoplus_{j=0}^{n-1}(\xi T_C)^{\otimes j},
 \qquad \rH^0(C,\mathcal A_b)=\bC,
\]
$\{x\}$ is $\mathcal A_b$-rigidifying in the sense of
\cite[Definition~3.6]{Heu25}. 
We denote by $\mathbf P_{\proet,\mathcal{A}_b}^x$ the $\{x\}$-rigidified pro-\'etale Picard functor of $\nu^*\mathcal{A}_b$ on $C$ \cite[Definition 3.12]{Heu25}. 
Put $R_x=\operatorname{Res}_{(\mathcal A_b|_x)/\bC}\mathbb G_m$.
By \cite[Definition~3.12 and Theorem~3.14(1)]{Heu25}, we have a short exact sequence:
\[
 1\longrightarrow R_x/\mathbb G_m
 \longrightarrow\mathbf P_{\proet,\mathcal{A}_b}^{x}
 \longrightarrow\mathbf P_{\proet,b}\longrightarrow1.
\]
The evaluation $\varepsilon_z$ induces $R_x\xrightarrow{\varepsilon_z^\times}\mathbb G_m$,
which splits the scalar subgroup $\Gm \subset R_x$. 
Thus we have a splitting $U_z:=\Ker(\varepsilon_z^\times)\xrightarrow{\sim}R_x/\mathbb G_m$. 
Then we have a natural isomorphism:
\[
\mathbf P_{\proet,\mathcal{A}_b}^{x}/U_z\xrightarrow{\sim}  \mathbf P_{\proet,b}^{z},\quad (\mathcal{L},\alpha:\mathcal{L}|_x\simeq \mathcal{A}_b)\mapsto (\mathcal{L}, \varepsilon_z^*\alpha). 
\]
The left hand side is also isomorphic to $\mathbf P_{\proet,b}$. 
By \cite[Theorem~2.4(4)]{Heu25}, we have
\[
 \Lie\mathbf P_{\proet,b}^{z}
 \simeq\rH^1(C_{\proet},\nu^*\mathcal A_b).
\]
This space also classifies additive $\nu^*\mathcal A_b$-torsors rigidified
at $z$.
Applying the functorial exponential \cite[Theorem~3.20]{Heu25}
(which extends \cite[Proposition~A.1.2]{Xu22}), we obtain an object in the neutral component of $\mathbf P_{\proet,b}^{z}$:
\[
 \Exp_{\mathbf P_{\proet,b}^{z}}
   \bigl([\mathcal P_{\mathcal C,b},\alpha_{\can}]\bigr)
 =[(\mathcal K_{\mathcal C,b}^{\Exp,z},\alpha_z)]
 \in (\mathbf P_{\proet,b}^{z})^{\circ}(\bC).
\]
The rigidified module is determined up to unique rigidified isomorphism since $\rH^0(C_{\proet},(\nu^*\mathcal A_b)^\times)=\bC^\times$ \cite[Lemma~3.11]{Heu25}.
Its image under the logarithmic Hodge--Tate map is
\[
 \mathrm{HTlog}
 (\mathcal K_{\mathcal C,b}^{\Exp,z})
 =\tau_b.
\]

Finally, by a similar argument of 
\cite[Proposition~A.2.4(i)]{Xu22}, one can show that $\mathcal K_{\mathcal C,b}^{\Exp,z}$ is indepenedent of the choice of $z$ up to a unique isomorphism. 
We denote by $\mathcal{K}_{\mathcal{C},b}^{\Exp}$, the resulting invertible $\nu^* \mathcal{A}_b$-module. 
Since $\mathrm{HTlog}\circ\Exp=\HT$, the module $\mathcal{K}_{\mathcal{C},b}^{\Exp}$ is an object of $\mathscr{H}_b$. 
\end{secnumber}

\begin{secnumber}
	Following \cite[Theorem~5.1(1)]{Heu25} and \cite[Theorem 6.2.5]{HX24}, the inverse of the $p$-adic Simpson correspondence defined by $(\mathcal{C},\Exp)$ is given, for a Higgs bundle $(E,\theta)$ above $b\in A(\bC)$, by
\begin{equation}\label{eq:inverse-spectral}
 \rS_{\mathcal C,\Exp}(E,\theta)
   =\nu^*E\otimes_{\nu^*\mathcal A_b}
       \mathcal K_{\mathcal C,b}^{\Exp}.
\end{equation}
This defines an exact tensor equivalence of categories $\rS_{\mathcal{C},\Exp}:\HB(C)\to \Bun_{\proet}(C)$. 

We denote its quasi-inverse, explicitly described in
\cite[Theorem~5.1(1)]{Heu25} and \cite[(6.2.12)]{HX24}, by
\begin{equation}\label{eq:HXExp-gen}
 \rH_{\mathcal C,\Exp}:\Bun_{\proet}(C)\xrightarrow{\sim}\HB(C). 
\end{equation}

Both sides admit a Hitchin map to the Hitchin base $A$ \cite{Heu-moduli,HX24}.
For each $b\in A(\bC)$, the fiber of above equivalence over $b$ is compatible with the natural action of $\mathscr{P}_b$ on both sides \cite[Theorem 6.2.5]{HX24}. 
\end{secnumber}

\begin{secnumber}\label{sss:common-spectral-algebra}
Let $f:C'\to C$ be a finite morphism of smooth projective $\bC$-curves, and $\mathcal C'$ a flat lifting
of $C'$ over $\rB_{\dR,2}^+$.
Write $\nu':C'_{\proet}\to C'_{\et}$ for the canonical morphism.
For $b\in A$ as above, recall \eqref{eq:C-fb} that 
\[
 \mathcal A_{f,b}=f^*\mathcal A_b,\qquad
 C'_{f,b}:=\Spec(\mathcal A_{f,b})\simeq C'\times_C C_b. 
\]
Let $b'=(b_i')_{i=1}^n$ be the image of $f^*b$ under
$f^*(\xi^{-i}\Omega_C^{\otimes i})\to\xi^{-i}\Omega_{C'}^{\otimes i}$.
The tangent map $df:\xi T_{C'}\to\xi f^*T_C$ induces natural maps
\begin{equation}\label{eq:common-spectral-map}
	\iota:\mathcal A_{b'}\longrightarrow\mathcal A_{f,b},
	\qquad q:C'_{f,b}\longrightarrow C'_{b'}.
\end{equation}
We compare the twisting objects over 
$\mathcal A_{f,b}$, and put
\[
 \mathcal K'=
 \mathcal K_{\mathcal C',b'}^{\Exp}
       \otimes_{\nu'^*\mathcal A_{b'}}\nu'^*\mathcal A_{f,b}.
\]
\end{secnumber}

\begin{prop}\label{p:twisting-comparison}
	\textnormal{(i)} With the above data, there is a canonical
	isomorphism of pro-\'etale line bundles on $C_{f,b}'$:
\begin{equation}\label{eq:twisting-line-comparison}
 f_{\proet}^*\mathcal K_{\mathcal C,b}^{\Exp}
       \otimes_{\nu'^*\mathcal A_{f,b}}(\mathcal K')^{-1}
 \xrightarrow{\sim}\nu'^*\mathcal L_{f,b}^{\Exp}.
\end{equation}

\textnormal{(ii)} The above isomorphism induces a natural isomorphism of functors $\Bun_{\proet}(C)\to \HB(C')$
\begin{equation}\label{eq:twisted-functoriality}
 \rH_{\mathcal C',\Exp}\circ f_{\proet}^*
 \xrightarrow{\sim}
 f^\circ_{\mathcal C',\mathcal C,\Exp}\circ\rH_{\mathcal C,\Exp}.
\end{equation}
\end{prop}

\begin{proof}
(i) Let $U=\Spa(R,R^+)$ be an affinoid perfectoid object of $C'_{\proet}$.
By \ref{sss:Higgs-Tate-torsor}, a section of $\mathcal{P}_{\mathcal{C}'}$ over $U$ is a $\rB_{\dR,2}^+$-algebra homomorphism $\psi:\mathscr{O}_{\mathcal{C}'}\to \BB_{\dR,2}^+(R)$ lifting $\mathscr{O}_{C'}\to R$, and a section of $f_{\proet}^*\mathcal{P}_{\mathcal{C}}$ over $U$ is such a homomorphism with source $f^{-1}\mathscr{O}_{\mathcal{C}}$, lifting the composition $f^{-1}\mathscr{O}_C\to \mathscr{O}_{C'}\to R$.
For a local section $\widetilde f$ of $\mathcal{L}_f$ and $\widetilde{f}^{\#}$ the associated homomorphism of algebras, 
the composition $\psi\circ\widetilde{f}^{\#}$ is therefore a section of $f_{\proet}^*\mathcal{P}_{\mathcal{C}}$.
Then $(\psi,\widetilde{f})\mapsto \psi\circ\widetilde{f}^{\#}$ induces a canonical isomorphism of $\nu'^*(\xi f^*T_C)$-torsors
\begin{equation}\label{eq:additive-lifting-comparison}
 f_{\proet}^*\mathcal P_{\mathcal C}
   \simeq(df)_*\mathcal P_{\mathcal C'}
   \wedge\nu'^*\mathcal{L}_f,
\end{equation}
where $\wedge$ denotes the contracted product of torsors.

Next, we push out \eqref{eq:additive-lifting-comparison} along $\beta:\xi f^*T_C\to\mathcal A_{f,b}$.
The left-hand side becomes $f_{\proet}^*\mathcal{P}_{\mathcal{C},b}$ \eqref{eq:spectral-HT-torsor}.
Since $\beta\circ df=\iota\circ\beta_{b'}$, the push-out of $(df)_*\mathcal{P}_{\mathcal{C}'}$ is the base change $\iota_*\mathcal{P}_{\mathcal{C}',b'}$ of the $\nu'^*\mathcal{A}_{b'}$-torsor $\mathcal{P}_{\mathcal{C}',b'}$ along $\iota$ \eqref{eq:common-spectral-map}. 
We obtain a canonical isomorphism of $\nu'^*\mathcal A_{f,b}$-torsors compatible with the canonical rigidifications
\begin{equation} \label{eq:compare-Afb-torsors}
	f_{\proet}^*\mathcal{P}_{\mathcal{C},b}\simeq \iota_*\mathcal{P}_{\mathcal{C}',b'}\wedge \nu'^*(\beta_*\mathcal{L}_f).
\end{equation}

We exponentiate the above isomorphism in a Picard group rigidified at one point.
Choose $x\in C(\bC)$ with $b_n(x)=0$ and $x'\in C'(\bC)$ above $x$.
Set
\[
 z=(x,0)\in C_b,\qquad z_f=(x',0)\in C'_{f,b},\qquad
 z'=q(z_f)=(x',0)\in C'_{b'}.
\]
The evaluation $\varepsilon_{z,f}=f^*(\varepsilon_z):\mathcal{A}_{f,b}|_{x'}\to \bC$ satisfies
\[
 \varepsilon_{z_f}\circ\beta|_{x'}=0,\qquad
 \varepsilon_{z_f}\circ\iota|_{x'}=\varepsilon_{z'}.
\]
All three $\nu'^*(\mathcal{A}_{f,b})$-torsors in \eqref{eq:compare-Afb-torsors} have canonical
rigidifications at $z_f$, and the isomorphism \eqref{eq:compare-Afb-torsors} is compatible with them. 

Let $G:=\mathbf P_{\proet,\mathcal{A}_{f,b}}^{z_f}$ be the pro-\'etale Picard functor of
$\nu'^*(\mathcal{A}_{f,b})$ on $C'_{\proet}$, rigidified at $z_f$. 
Since $\rH^0(C',\mathcal{A}_{f,b})=\bC$, the argument of \ref{sss:twisting-object}
applies to $\mathcal{A}_{f,b}$ and shows that $G$ is representable by a rigid group. 
Pullback, extension of coefficients, and passage from the
\'etale site to the pro-\'etale site give homomorphisms
\[
 \mathbf P_{\proet,b}^{z}\longrightarrow G,\qquad
 \mathbf P_{\proet,b'}^{z'}\longrightarrow G,\qquad
 \mathbf P_{\et,f,b}^{z_f}\longrightarrow G.
\]
Via the their differentials, we may interpret \eqref{eq:compare-Afb-torsors} as an identity in $\Lie G$ compatible with rigidification:
\[
[f_{\proet}^*\mathcal{P}_{\mathcal{C},b}]=[\iota_*\mathcal{P}_{\mathcal{C}',b'}]+ [\nu'^*(\beta_*\mathcal{L}_f)].
\]
By applying the exponential $\Exp_G$, we obtain an isomorphism \eqref{eq:twisting-line-comparison} compatible with canonical rigidifications. 
This isomorphism is canonical and more precisely independent of $x$ by the argument
\cite[Proposition~A.2.4(i)]{Xu22}.

(ii) 
We apply the tensor product of \eqref{eq:twisting-line-comparison} with
$\nu'^*f^*E$. By \eqref{eq:inverse-spectral}, we deduce a canonical isomorphism
\[
 f_{\proet}^*\rS_{\mathcal C,\Exp}(E,\theta)
 \simeq
 \rS_{\mathcal C',\Exp}
     \bigl(f^\circ_{\mathcal C',\mathcal C,\Exp}(E,\theta)\bigr).
\]
This finishes the proof of \eqref{eq:twisted-functoriality}.
\end{proof}

\subsection{Comparison of torsors}

\begin{secnumber}
In this subsection, we consider rank two nilpotent Higgs bundles on $C$.
They lie over the zero point $b=0$ of the Hitchin base. 
Let 
\[\pi:C_0=\Spec_{\mathscr{O}_C}(\mathscr{O}_C\oplus \xi T_C)\to C\]
be the degree two spectral curve of $C$ over $b=0\in A(\bC)$, where $\xi T_C$ is an ideal of square zero. 
The natural surjection $\mathscr{O}_{C_0}\to \mathscr{O}_{C}$ defines a section $\iota:C\to C_0$ of $\pi$, identifying $C$ with the reduced subscheme of $C_0$. 
This splitting yields a canonical decomposition of the group of isomorphism classes of line bundles on $C_0$:
\begin{equation}
	\rH^1(C_0,\mathscr{O}_{C_0}^{\times})\simeq \rH^1(C,\mathscr{O}_C^{\times})\times \rH^1(C,1+ \xi T_C). 
	\label{eq:decomposition-Pic-C0}
\end{equation}
The first factor consists of pullbacks along $\pi$, and the second of line bundles trivial on $C$.
\end{secnumber}

\begin{secnumber}
We consider the following torsors under the cohomology group $\rH^1(C,\xi T_C)$: 

\begin{enumerate}
	\item The torsor of isomorphism classes of flat liftings of $C$ over $\rB_{\dR,2}^+$.
To fix convention, for two liftings $\mathcal C_1,\mathcal C_2$
we write $\mathcal L_{\mathcal C_1,\mathcal C_2}$ for the torsor of liftings
$\mathcal C_1\to\mathcal C_2$ of $\id_C$ on $C_{\et}$, namely the torsor
$\mathcal L_{\id}$ of \ref{sss:review-LExp} with source $\mathcal C_1$,
target $\mathcal C_2$, and $f=\id_C$.
Our convention is
\[
 \mathcal C_2-\mathcal C_1
       =[\mathcal L_{\mathcal C_1,\mathcal C_2}]
       \in\rH^1(C,\xi T_C).
\]

	\item The torsor of equivalence classes of extensions of pro-\'etale vector bundles on $C_{\proet}$: 
		\begin{equation} \label{eq:ext-Faltings}
		0\to \xi \hatO_C \to \mathcal{E} \to \nu^*\Omega_C \to 0,
		\end{equation}
		whose class $[\mathcal{E}]\in \rH^1(C_{\proet},\nu^*(\xi T_C))$ maps to $1\in \rH^0(C,\mathscr{O}_C)$ under the Hodge--Tate map \eqref{eq:HT-sequence-W}
		\[
		0\to \rH^1(C, \xi T_C) \to \rH^1(C_{\proet}, \nu^*(\xi T_C)) \to \rH^0(C, \mathscr{O}_C) \to 0.
		\]

	\item The torsor of isomorphism classes of line bundles on $C_0$ whose restriction to $C$ is trivial. 

	\item The torsor of equivalence classes of extensions of Higgs bundles on $C$:
		\begin{equation} \label{eq:ext-Higgs}
		0\to (\xi\mathscr{O}_C,0) \to (E,\Theta) \to (\Omega_C,0) \to 0,
	\end{equation}
	such that the Higgs field $\Theta:E\to E\otimes_{\mathscr{O}_C}\xi^{-1}\Omega_C$ is induced by $\id_{\Omega_C}: \Omega_C \to (\xi \mathscr{O}_C)\otimes_{\mathscr{O}_C}\xi^{-1}\Omega_C$. 
\end{enumerate}
\end{secnumber}
\begin{prop} \label{p:Torsors-equiv}
	\textnormal{(i)} There exist canonical isomorphisms between torsors (1) $\Leftrightarrow$ (2) and (3) $\Leftrightarrow$ (4). 

	\textnormal{(ii)} 
	A lift $\mathcal{C}$ of $C$ induces isomorphisms of torsors (1) $\Leftrightarrow$ (3) and (2) $\Leftrightarrow$ (4). 
\end{prop}

\begin{proof}
(i)
(1) $\Leftrightarrow$ (2).  A flat lift $\mathcal{C}$ of $C$ over $\rB_{\dR,2}^+$ defines the Higgs--Tate torsor $\mathcal{P}_{\mathcal{C}}$ under $\nu^*(\xi T_C)$ on $C_{\proet}$ \ref{sss:Higgs-Tate-torsor}. 
Its associated class $[\mathcal{P}_{\mathcal{C}}]$ in $\rH^1(C_{\proet},\nu^*(\xi T_C))$ has the Hodge--Tate image $1\in \rH^0(C,\mathscr{O}_C)$ \cite[Proposition 2.15]{Heu25}. 
This defines an extension as in \eqref{eq:ext-Faltings}. 

To check equivariance, 
for two liftings $\mathcal C',\mathcal C$, 
a local section $f:\mathcal C'\to\mathcal C$ of $\mathcal{L}_{\mathcal{C}',\mathcal{C}}$ induces $f^{\#}:\mathscr O_{\mathcal C}\to\mathscr O_{\mathcal C'}$.
The composition $\varphi\mapsto\varphi\circ f^{\#}$ identifies the contracted sum of $P_{\mathcal C'}$ and $\nu^*\mathcal{L}_{\mathcal{C}',\mathcal{C}}$ with $P_{\mathcal C}$ \eqref{eq:additive-lifting-comparison}.
Thus, we have
\begin{equation}\label{eq:variation-dual-Faltings}
	[\mathcal{P}_{\mathcal C}]-[\mathcal{P}_{\mathcal C'}]
	=[\nu^*\mathcal{L}_{\id}].
\end{equation}
Hence $\mathcal{C}\mapsto \mathcal{P}_{\mathcal{C}}$ is equivariant and then a bijection. 

(3) $\Leftrightarrow$ (4). We show that the functor $\mathcal{L}\mapsto \pi_*(\mathcal{L})\otimes_{\mathscr{O}_C} \Omega_C$ defines an $\rH^1(C,\xi T_C)$-equivariant bijection between (3) and (4). 
Given a line bundle $\mathcal{L}$ on $C_0$ equipped with a trivialization $\iota^*\mathcal{L}\simeq \mathscr{O}_C$, tensoring the short exact sequence $0\to \xi T_C\to \mathscr{O}_{C_0}\to \mathscr{O}_C\to 0$ with $\mathcal{L}$ and applying $\pi_*$ yields a short exact sequence of $\mathscr{O}_C$-modules
\[
	0\to \pi_*(\mathcal{L}\otimes_{\mathscr{O}_{C_0}} \xi T_C) \to \pi_*(\mathcal{L}) \to \pi_*(\mathcal{L}\otimes_{\mathscr{O}_{C_0}}\mathscr{O}_C)(\simeq \mathscr{O}_C) \to 0.
\]
The tautological $\mathscr{O}_{C_0}$-module structure on $\mathcal{L}$ equips $\pi_*(\mathcal{L})$ with a Higgs field $\theta:\pi_*(\mathcal{L})\to \pi_*(\mathcal{L})\otimes_{\mathscr{O}_C}\xi^{-1}\Omega_C$ that factors through an isomorphism $\mathscr{O}_C \xrightarrow{\sim} \pi_*(\mathcal{L}\otimes_{\mathscr{O}_{C_0}} \xi T_C) \otimes_{\mathscr{O}_C} \xi^{-1}\Omega_C$. 
Twisting by $\Omega_C$, we see that $\pi_*(\mathcal{L})\otimes_{\mathscr{O}_C}\Omega_C$ defines an object of (4). 
Conversely, given an extension $(E,\theta)$ as in (4), the Higgs field $\theta$ endows $E$ with the structure of an invertible $\mathscr{O}_{C_0}$-module such that $\iota^*(E\otimes_{\mathscr{O}_C}T_C)\simeq \mathscr{O}_C$. Thus $E\otimes_{\mathscr{O}_C} T_C$ defines a quasi-inverse from (4) to (3).

(ii) (1) $\Leftrightarrow$ (3). Fix a flat lifting $\mathcal{C}$ of $C$ over $\rB_{\dR,2}^+$. For any flat lifting $\mathcal{C}'$ of $C$ over $\rB_{\dR,2}^+$, the sheaf of flat liftings of $\id_C$ to morphisms $\mathcal{C}\to \mathcal{C}'$ forms a $\xi T_C$-torsor $\mathcal{L}_{\mathcal{C},\mathcal{C}'}$ over $C$. 
Pushing out $\mathcal{L}_{\mathcal{C},\mathcal{C}'}$ along $\xi T_C \to \mathscr{O}_{C_0}^{\times},\quad u\mapsto 1+u$, we obtain an object of (3). 
This assignment is $\rH^1(C,\xi T_C)$-equivariant and is therefore an isomorphism of torsors. 

(2) $\Leftrightarrow$ (4). 
Given a flat lifting $\mathcal{C}$,
if an object of (4) has underlying extension class
$e\in\rH^1(C,\xi T_C)$, then \eqref{eq:inverse-spectral} shows that its image under
$\rS_{\mathcal C}$ has extension class
$\nu^*e+[\mathcal{P}_{\mathcal C}]$. 
Hence, the functor $\rH_{\mathcal{C},\Exp}$ induces an $\rH^1(C,\xi T_C)$-equivariant bijection between (2) and (4). 
\end{proof}

\begin{rem}
The Faltings extension $\mathcal{F}_{\mathcal{C}}$ is the sheaf on $C_{\proet}$ of affine functions on $\mathcal{P}_{\mathcal{C}}$ \cite[II.10.13]{AGT}:
\begin{equation} \label{eq:affine-function}
\mathcal{F}_{\mathcal{C}}:=\textnormal{Aff}_{\widehat{\mathscr{O}}_C}(\mathcal{P}_{\mathcal{C}},\widehat{\mathscr{O}}_C). 
\end{equation}

It fits into a short exact sequence:
\begin{equation}\label{eq:Faltings-seq}
0\to \widehat{\mathscr{O}}_C\to \mathcal{F}_{\mathcal{C}}\to \nu^*(\xi^{-1} \Omega_C)\to 0
\end{equation}
where the first map is the constant map and the second map is the linear part. 
By \eqref{eq:affine-function}, we deduce that dual of Faltings extension $\mathcal{F}_{\mathcal{C}}$ is isomorphic to $\mathcal{P}_{\mathcal{C}}$ as $\nu^*(\xi T_C)$-torsor. We refer to \cite[Remark 2.17]{Heu25} for another argument. 
In particular, $\HT(\mathcal{F}_{\mathcal{C}})=-1$. 
\end{rem}
\begin{rem}
	A characteristic $p$ analogue of \Cref{p:Torsors-equiv} can be found in \cite[Theorem 4.5]{OV07} and \cite[\S~3.5]{CZ15}. 
\end{rem}

\section{Extensions of line bundles}
\label{s:extensions}
Let $C$ be a smooth projective curve over an algebraically closed field $F$ of genus $\ge 2$ and $L$ a line bundle of degree $d \ge 1$. The goal of this section is the following. 

\begin{prop} \label{p:ext-L}
	\textnormal{(i)} There exists extensions $E$ in 
\begin{equation}
	0\to L^{-1} \to E\to L\to 0,
	\label{eq:extension}
\end{equation}
such that $E$ is semi-stable and non-stable. 
In particular, $E$ is an extension of line bundles of degree zero. 

\textnormal{(ii)} When $d\ge 2$, the projective extension classes whose middle terms are semi-stable and non-stable form a locally closed subvariety in $\mathbb{P}(\rH^1(C,L^{\otimes -2})^{\vee})$ of dimension $2d-1$.
\end{prop} 

\subsection{Case \texorpdfstring{$d=1$}{d=1}}

In this case, \Cref{p:ext-L} follows from the following:

\begin{prop}[\cite{NR}] \label{p:ext-squareroot} Suppose $d=1$. 
	
	\textnormal{(i)} Any extension $E$ in \eqref{eq:extension} is either semi-stable or $E=L\oplus L^{-1}$. 

	\textnormal{(ii)} There exists an extension $E$ as above such that the $E$ is semi-stable and non-stable. 
\end{prop}

\begin{proof}
	Assertion (i) follows from \cite[Lemma 5.1]{NR}. To show assertion (ii), we recall the following result. 

\begin{lemma}[\cite{NR} Lemma 5.1, 5.2]
	A non-trivial extension $E$ in \eqref{eq:extension} is non-stable if and only if there is a closed point $x\in C$ such that its extension class $\delta(E)$ lies in the kernel of the map
	\[
	\rH^1(X,L^{\otimes -2}) \to \rH^1(X,L^{\otimes -2}\otimes L_x),
	\]
	where $L_x$ is the line bundle associated to the effective divisor $x$. In this case, $E$ is $S$-equivalent to $L\otimes L_x^{-1} \oplus L^{-1}\otimes L_x$. 
\end{lemma}

	(ii) For a point $x\in C$, we consider the long exact sequence:
	\[
	0\to \rH^0(C, L^{\otimes -2}) \to \rH^0(C,L^{\otimes -2}\otimes L_x) \to \rH^0(C,F(x)) \to \rH^1(C,L^{\otimes -2}) \to \rH^1(C,L^{\otimes -2}\otimes L_x),
	\]
	where $F(x)=(L^{\otimes-2}\otimes L_x)|_x$ is a skyscraper sheaf, the first and the second terms are zero, and the third is one-dimensional. Hence, we may choose an extension $E$ whose nonzero class $\delta(E)$ lies in the kernel of the last map. Then assertion (ii) follows from the above lemma. 
\end{proof}

\subsection{Secant varieties associated to a non-degenerate embedding of a curve,
after Lange \texorpdfstring{\cite{Lange}}{}}

Let $\varphi:C\to \mathbb{P}^N$ be a closed embedding of $C$ into a projective space of dimension $N$ over $F$. 
Let $n$ be an integer with $1\le n\le N$. 
An \textit{$n$-secant divisor} of $C$ is an effective divisor of degree $n$ of $X$, whose linear span in $\mathbb{P}^N$ is of dimension $n-1$. 
The linear span itself is called an \textit{$n$-secant} of $C$. 
The \textit{$n$-secant variety} $\Sec_n(C)$ is the closure of the union of all $n$-secants of $C$ in $\mathbb{P}^N$. 
We have a sequence of embeddings:
\[
	C=\Sec_1(C) \subset \Sec_2(C) \subset\cdots \subset \Sec_n(C) \cdots
\]
In general, we have $\dim \Sec_n(C)\le 2n-1$. 

\begin{theorem}[\cite{Lange} \S~1 Theorem]
	Suppose $\varphi$ is non-degenerate, i.e. $C$ is not contained in any hyperplane. 
	Then we have $\dim \Sec_n(C)=\min\{ 2n-1,N\}$. 
\end{theorem}

\subsection{Case \texorpdfstring{$d\ge 2$}{d greater than or equal to 2}}
\begin{secnumber}
	We may rewrite the extension \eqref{eq:extension} as
\[
0\to \mathscr{O}_C \to E\otimes L \to L^{\otimes 2}\to 0. 
\]

Since $\deg L^{\otimes 2}\otimes \Omega_C=2g-2+2d \ge 2g+2$, the line bundle $L^{\otimes 2}\otimes \Omega_C$ is very ample and the map 
\begin{equation} \label{eq:embedding}
\varphi: C\to \mathbb{P}=\mathbb{P}(\rH^0 (C, L^{\otimes 2}\otimes \Omega_C)), \quad x \mapsto \{ \sigma\in \rH^0(C,L^{\otimes 2}\otimes \Omega_C)| ~ \sigma(x)=0 \}
\end{equation}
associated to the complete linear system $|L^{\otimes 2}\otimes \Omega_C|$ is a closed embedding. 
Moreover, $C$ is non-degenerately embedded into $\mathbb{P}$.
Note that $\mathbb{P}$ has dimension
\[N=g-2+2d.\]

	For any non-trivial extension \eqref{eq:extension}, we consider the extension class $\delta(E)\in \rH^1(C,L^{\otimes -2})\simeq \rH^0(C,L^{\otimes 2}\otimes \Omega_C)^{\vee}$ as a point of $\mathbb{P}$. 
\end{secnumber}
\begin{secnumber}
	Let $E$ be a vector bundle of rank $2$ on $C$, we define the invariant $s(E)$ by 
	\[
	s(E)=\deg(E)-2 \max \deg(M),
	\]
	where the maximum being taken over all line subbundles $M$ of $E$.  
	By definition, $E$ is stable (resp. semi-stable) if and only if $s(E)>0$ (resp. $s(E)\ge 0$). 

	Let $E$ be a non-trivial extension in \eqref{eq:extension}. Then $s(E)$ is an even number and we have $2-2d\le s(E) \le 2d$. 
	The geometric meaning of $s(E)$ is related to the secant varieties of $\varphi$ \eqref{eq:embedding}. 
\end{secnumber}

\begin{prop}[\cite{LN} Proposition 1.1] \label{p:Lange-Narasimhan}
	Let $s$ be an even integer such that $4-2d\le s \le 2d$. 
	The following conditions are equivalent:

	\textnormal{(i)} $s(E) \ge s$. 

	\textnormal{(ii)} $\delta(E)\notin \Sec_{d-1+\frac{s}{2}}(C)$.
\end{prop}

\begin{secnumber}
	\textit{Proof of \Cref{p:ext-L}}. 
	When $s=0$ (resp. $s=2$) and $n=d-1$ (resp. $n=d$), we have $\dim \Sec_{d-1}(C)=\min\{ 2n-1,N=2d+g-2\}=2d-3$ (resp. $\dim \Sec_{d}(C)=2d-1$). 
	
	For any extension $E$ with class $\delta(E)\in \Sec_d(C)$ but $\delta(E)\notin \Sec_{d-1}(C)$, we have $s(E)=0$ by \Cref{p:Lange-Narasimhan}.
	Such a vector bundle $E$ is semi-stable and non-stable. 
	And the dimension of the locally closed subvariety $\Sec_d(C)-\Sec_{d-1}(C)$ is $2d-1$. \hfill \qed 
\end{secnumber}

\section{Proof of \Cref{t:intro-nilpotent}}
\label{s:main-result}

\subsection{Preliminaries}\label{ss:section4-preliminaries}
Let $L$ be a line bundle of degree $d \ge 0$ on $C$. 
The group of equivalence classes of extensions of Higgs bundles:
\begin{equation} \label{eq:ext-Higgs-L}
0\to (L^{-1},0) \to (E,\Phi) \to (L,0)\to 0
\end{equation}
is isomorphic to 
\[
\rH^1(C,L^{-2})\oplus \rH^0(C, L^{-2}\otimes \Omega_C). 
\]
The first component measures the extension of underlying bundles and the second classifies Higgs fields. 

Let $\eta:L^{\otimes 2}\to \Omega_C$ be a section of $\rH^0(C, L^{-2}\otimes \Omega_C)$. 
For an extension class $\gamma\in \rH^1(C,L^{-2})$, 
we denote by $E(\gamma,\eta)$ the Higgs bundle in \eqref{eq:ext-Higgs-L} corresponding to $(\gamma,\eta)$:
\[
0\to (L^{-1},0) \to \left(E_{\gamma},\begin{pmatrix}0& \xi^{-1} \eta \\0&0\end{pmatrix}\right) \to (L,0)\to 0. 
\]
With the above notation, the uniformizing Higgs bundle corresponds to $E(0,\theta)$. 

For $a\in\bC^\times$, pushing out the extension defining $E_\gamma$ by multiplication by $a$ on $L^{-1}$ gives an extension of class $a\gamma$ and an isomorphism $E_{\gamma}\xrightarrow{\sim} E_{a\gamma}$.
This carries the Higgs field $\xi^{-1}\eta$ to $\xi^{-1}a\eta$ and induces an isomorphism:
\begin{equation} \label{l:section4-scaling}
E(\gamma,\eta)\xrightarrow{\sim} E(a\gamma,a\eta). 
\end{equation}
In particular, all $E(0,a\eta)$ with $a\ne0$ are isomorphic. 
If we fix $\eta\neq 0$, the following natural map is injective:
\[
\Ext^1(L,L^{-1})\to \{\textnormal{isomorphism classes of Higgs bundles}\},\qquad \gamma\mapsto E(\gamma,\eta).
\]

The dual $\eta^{\vee}:T_C\to L^{-2}$ of $\eta$ is injective and its quotient has finite support. Then we have a surjection:
\begin{equation}
	\mu_{\eta}:=\rH^1(\eta^{\vee}): \rH^1(C,T_C)\twoheadrightarrow \rH^1(C,L^{-2}). 
	\label{eq:surjection-mu}
\end{equation}
For every $a\in \bC^{\times}$, we have $\mu_{a\eta}=a\mu_{\eta}$. 

\begin{prop}
\label{p:twisted-action-rank-two}
Let $\mathscr{L}$ be a line bundle on the square-zero spectral
curve $C_0$ with trivial
restriction to $C$. 
Write $[\mathscr{L}]\in\rH^1(C,1+\xi T_C)$ and
$u=\xi^{-1}[\mathscr{L}]\in\rH^1(C,T_C)$.
Then we have
\begin{equation}\label{eq:section4-twisted-action}
 E(\gamma,\eta)
       \otimes_{\mathscr{O}_{C_0}}\mathscr{L}
 \simeq
 E(\gamma+\mu_{\eta}(u),\eta).
\end{equation}
\end{prop}

\begin{rem}
	The dependence of the $p$-adic Simpson correspondence on the choice of lifting has recently been investigated by Pan--Yu \cite{PY}. 
	In our setting, \Cref{p:Torsors-equiv} and \Cref{p:twisted-action-rank-two} give an explicit description of the dependence on the lifting for rank-two nilpotent Higgs bundles. 
\end{rem}

\begin{proof}
	Choose local splittings $\varphi_i:L^{-1}\oplus L\xrightarrow{\sim}E_\gamma$ and local generators $s_i$ of $\mathscr L$ on an open subset $U_i$ of $C$. Then we can represent the class $\gamma$, $u$ as \v{C}ech cocycles $(\gamma_{ij}\in L^{-2}(U_i))$, $(u_{ij}\in T_{C}(U_i))$:
\[
\varphi_i^{-1}\varphi_j=\begin{pmatrix}
	 1& \gamma_{ij}\\
  0&1
 \end{pmatrix}\qquad 
 s_j=(1+\xi u_{ij})s_i. 
\]
%The spectral action of the element $\xiu_{ij}$ on $\gamma_{ij}$ by $c\eta^\vee(u_{ij}))$. 
The spectral action of $\mathscr{L}$ on $E(\gamma,c)$ is given by the following formula on 
the transition maps:
\[
 \begin{pmatrix}
	 1& \xi u_{ij}\\
  0&1
 \end{pmatrix}\cdot 
 \begin{pmatrix}
	 1& \gamma_{ij}\\
  0&1
 \end{pmatrix}
 =
 \begin{pmatrix}
	 1& \gamma_{ij}+\eta^{\vee}(u_{ij})\\
  0&1
 \end{pmatrix}.
\]
Hence, the underlying bundle in the right hand side has the extension class $\gamma+\mu(u)$. Moreover, this action does change the Higgs field locally. 
This proves \eqref{eq:section4-twisted-action}.
\end{proof}

\begin{lemma}\label{l:extensions-degree-zero}
Fix a lifting $\mathcal C$ and an exponential $\Exp$.
Let $M_1,M_2$ be line bundles of degree zero on $C$, and let $W_i$ be the
continuous representations corresponding to $(M_i,0)$ \cite[Theorem 1.1]{Heu24}.
The natural map
\[
 \Ext^1_{\Rep(\pi_1(C,\overline x),\bC)}(W_2,W_1)
 \longrightarrow
 \Ext^1_{\HB(C)}((M_2,0),(M_1,0))
\]
is an isomorphism. In particular, every Higgs extension of $(M_2,0)$ by
$(M_1,0)$ lies in the essential image of $\bH_{\mathcal C,\Exp}$.
\end{lemma}

\begin{proof}
We set $W=\Hom_{\bC}(W_2,W_1)$, and
$M=M_1\otimes M_2^{-1}$. The pro-\'etale vector bundle associated with $W$
is $\nu^*M$. Let $q:C^{\univ}\to C$ be the universal
pro-finite-\'etale cover. 
By the primitive comparison theorem, we have a natural isomorphism
\[
 \rH^1_{\cont}(\pi_1(C,\overline x),W)
 \xrightarrow{\sim}\rH^1_v(C,\nu^*M). 
\]
The two sides classify extensions of $W_2$ by $W_1$ and extensions of
$\nu^*M_2$ by $\nu^*M_1$, respectively. By the cohomology comparison \cite[Theorem~5.5]{Heu25}, we have
\[
 \rH^1_{\cont}(\pi_1(C,\overline x),W)
 \simeq \rH^1(C,M)\oplus
 \rH^0(C,M\otimes\xi^{-1}\Omega_C).
\]
The right hand side classifies the extensions of $(M_2,0)$ by $(M_1,0)$ and the lemma follows. 
\end{proof}

\subsection{Openness of the representation locus}\label{ss:good-case}

\begin{secnumber}
	We recall the openness of the representations locus in \cite{HX24}.
	Let $\bfHig_{2}$ and $\bfBun_{2,v}$ denote the $v$-sheaves obtained by sheafifying the functors of isomorphism classes of rank-two Higgs bundles and rank-two vector $v$-bundles,
respectively.  Their $\bC$-points are sets of isomorphism classes of rank-two Higgs bundles (resp. vector $v$-bundles) on $C$. 

We recall the topology on the sets of $\bC$-points in \cite[Definition~7.4.1]{HX24}. 
For a profinite set $S$, we set
\[
 T_S=\Spa\bigl(\operatorname{Cts}(S,\bC),
                    \operatorname{Cts}(S,\oo)\bigr).
\]
For either of the above $v$-sheaves $F$, 
we endow $F(\bC)$ the finest topology such that for every morphism $T_S\to F$, the associated map on $\bC$-points $S\to F(\bC)$ is continuous. 
To avoid set-theoretic issues, we fix a cardinal $\kappa$ and work with $\kappa$-small pro-finite sets $S$. 

If we fix $\mathcal{C}$ and an exponential $\Exp$ for $\bC$, then moduli version of the $p$-adic Simpson correspondence \cite{HX24} induces a homeomorphism \cite[Theorem 7.4.2]{HX24}:
\begin{equation}
	S_{\mathcal{C},\Exp}: \bfHig_{2}(\bC)\xrightarrow{\sim} \bfBun_{2,v}(\bC). 
	\label{eq:top-hom}
\end{equation}
By \cite[Theorem~9.2.1 and Remark~9.2.2]{HX24}, the
$v$-sheaf of conjugacy classes of continuous representations is
an open sub-$v$-sheaf:
\[\bfRep_{2}\to \bfBun_{2,v}.\]
\end{secnumber}

\begin{lemma}
	We set $V:=\rH^1(C,L^{-2}), W:=\rH^0(C,L^{-2}\otimes\Omega_C)$.  
The natural map
\[
 V\times W\longrightarrow\bfHig_{2}(\bC),
 \qquad (\gamma,\eta)\longmapsto E(\gamma,\eta),
\]
is continuous for the usual $p$-adic topology on a
finite-dimensional $\bC$-vector space.
\end{lemma}

\begin{proof}
Let $S=\mathbb N\cup\{\infty\}$ be the profinite space in which
the neighborhoods of $\infty$ contain all but finitely many
integers.  A convergent sequence in $V\times W$, together with its limit, defines a continuous map 
\[S\to V\times W,\quad n\mapsto (\gamma_n,\eta_n), \infty \mapsto (\gamma_{\infty},\eta_{\infty}).\]
The family of Higgs bundles $E(\gamma_n,\eta_n)$ converges to $E(\gamma_{\infty},\eta_{\infty})$ and defines a Higgs bundle on $C\times_{\bC} T_S$ and therefore a section $z\in \bfHig_{2}(T_S)$. 
In particular, the composition $S\to V\times W\to \bfHig_{2}(\bC)$ is continuous. 
Since $V\times W$ is metrizable, the sequential continuity implies the continuity of $V\times W\to \bfHig_{2}(\bC)$. 
\end{proof}

\begin{proof}[Proof of \Cref{t:intro-nilpotent}]
	Recall that any line bundle of degree zero, equipped with the zero Higgs field, lies in the essential image of $\bH_{\mathcal{C}}$ \cite[Theorem~1.1]{Heu24}, \cite[Theorem~12]{DW05}.
	Let $(E,\Theta)$ be a rank two semi-stable Higgs bundle of degree zero with non-zero nilpotent Higgs field over $C$. 
	To prove the theorem, after tensor with a line bundle, we may assume that $\det(E,\Theta)\simeq (\mathscr{O}_{C},0)$. 
	
	Set $L^{-1}=\Ker(\Theta)$. Since $(E,\Theta)$ is semi-stable, $L$ is a line bundle of degree $d\ge 0$, and $(E,\Theta)$ fits into an extension \eqref{eq:ext-Higgs-L} defined by a non-zero section $\eta:L^{\otimes 2}\to \Omega_C$. 
	Thus $(E,\Theta)\simeq E(\gamma,\eta)$ for some extension class $\gamma\in \rH^1(C,L^{-2})$.

	We take a lift $\mathcal{C}$ of $C$ over $\rB_{\dR,2}^+$ and an exponential $\Exp$. 
By \Cref{p:ext-L}, there is a class
$\gamma_0\in\rH^1(C,L^{-2})$ for which $E_{\gamma_0}$ is an extension
of two line bundles of degree zero. Hence
$E(\gamma_0,0)$ lies in
$S_{\mathcal C,\Exp}^{-1}(\bfRep_2(\bC))$
by \Cref{l:extensions-degree-zero}.
This subset is open, and the map
$a\mapsto E(\gamma_0,a\eta)$ is continuous.
Consequently, for some $a\in\bC^\times$ sufficiently close to zero,
$E(\gamma_0,a\eta)$ belongs to this subset.

	Since $\mu_{a\eta}:\rH^1(C,T_C)\twoheadrightarrow \rH^1(C,L^{-2})$ is surjective by \eqref{eq:surjection-mu}, there exists $u\in \rH^1(C,T_C)$ such that
	\[
	\mu_{a\eta}(u)=a\gamma-\gamma_0.
	\]
	By \Cref{p:Torsors-equiv}, there exists another flat lifting $\mathcal{C}'$ of $C$ over $\rB_{\dR,2}^+$ such that the torsor $\mathcal{L}_{\mathcal{C}',\mathcal{C}}$ has class $\xi^{-1}[\mathcal{L}_{\mathcal{C}',\mathcal{C}}]=u$.
	Applying \Cref{p:twisting-comparison}(ii) and \Cref{p:twisted-action-rank-two}, together with the scaling isomorphism \eqref{l:section4-scaling}, we deduce that
	\[
	E(\gamma,\eta)\simeq E(a\gamma,a\eta)\simeq E(\gamma_0,a\eta)\otimes_{\mathscr{O}_{C_0}}\mathcal{L}_{\mathcal{C}',\mathcal{C}}^{\Exp}
	\]
	lies in $S_{\mathcal{C}',\Exp}^{-1}(\bfRep_{2}(\bC))$.
This finishes the proof. 
\end{proof}

\begin{proof}[Proof of corollary 1.2.2 \Cref{c:kostant}] 
For $n\ge 2$, we have an isomorphism:
\[
 (E_n,\Theta_n)\xrightarrow{\sim}\Sym^{n-1}(E_2,\Theta_2).
\]
Then the corollary follows from \Cref{t:intro-nilpotent} and the fact that $\bH_{\mathcal{C},\Exp}$ is a tensor equivalence. 
\end{proof}

\subsection{Comments} \label{ss:comments}
\begin{secnumber} \label{sss:alternative-small}
	The proof of \Cref{t:intro-nilpotent} in \S~\ref{ss:good-case} modifies the flat lifting of $C$ over $\rB_{\dR,2}^+$ in order to adjust the extension class $\gamma$, and therefore does not prove \Cref{t:intro-nilpotent} for every lifting.

	When $C$ descends to $K$ and admits a semi-stable $\OK$-model $X$, the parallel transport construction of \cite{Xu22} suggests an alternative approach to \Cref{t:intro-nilpotent} via the small $p$-adic Simpson correspondence over $X$ \cite{AGT,AGrel}.
	This approach provides a finer control on the lifting. 
	The idea is to start from a semi-stable, non-stable bundle $E_{\gamma_0}$, which admits a strongly semi-stable reduction and hence corresponds to a representation by \cite{DW05}. Then we equip it with a sufficiently small nilpotent Higgs field $a\eta$ so that it falls into the image of $\mathbb{H}_{\mathcal{C},\Exp}$ (after passing to a suitable finite cover, as in \cite{Xu22}). 
	However, recovering $E(\gamma,\eta)$ from $E(\gamma_0,a\eta)$ still requires varying the lifting of $C$ over $\rB_{\dR,2}^+$ to adjust the extension class. 
	Thus, both arguments share the same limitation: neither proves \Cref{t:intro-nilpotent} for every lifting. 

	On the other hand, for Shimura and modular curves, the uniformizing Higgs bundle (or its logarithmic variant) is known to correspond to a local system under the canonical lifting, as we recall in \ref{sss:canonical-modular}. This follows from the $p$-adic Riemann--Hilbert correspondence \cite{LZ17,DLLZ23}. 
	However, as we illustrate in \ref{sss:modular-tower-obstruction}, even in the modular curve case one cannot deduce this correspondence directly from the small $p$-adic Simpson correspondence along the natural tower of canonical liftings. 
\end{secnumber}

\begin{secnumber}\label{sss:canonical-modular}
Put $K_0=\operatorname{Frac}W(k)$.
Let $N\ge3$ be prime to $p$, fix a primitive $N$-th root of unity
$\zeta_N\in K_0$, and let $X=X_N$ over $K_0$ be the connected compactified modular
curve with Weil pairing $\zeta_N$. Write $j:U=X-D\hookrightarrow X$
for the complement of its cusp divisor.
If $h:\mathcal E\to U$ is the universal elliptic curve, we set
\[
 \mathbb L=j_{\mathrm{k\acute et},*}R^1h_{\et,*}\Qp,
 \qquad
 \mathbb L^\circ=j_{\mathrm{k\acute et},*}R^1h_{\et,*}\Zp.
\]
the Kummer–étale local system on $X$ \cite[\S~4.1]{DLLZ23}. 

Applying logarithmic $p$-adic Riemann--Hilbert over $K_0$ to the de Rham local system $\mathbb L$, we obtain the logarithmic Gauss--Manin connection $(\mathcal H,\nabla,\mathrm{Fil}^{\bullet})$, equipped with a Hodge filtration \cite[Theorem~5.3.1]{DLLZ23}:
\[
 0\longrightarrow\omega=\mathrm{Fil}^1\mathcal H
 \longrightarrow\mathcal H\longrightarrow\omega^{-1}\longrightarrow0.
\]
Griffiths transversality induces a Higgs field on the associated graded:
\[
 (E,\vartheta)=
 \bigl(\omega_{\bC}(-1)\oplus\omega_{\bC}^{-1},\gr\nabla\bigr).
\]
This Higgs field is induced by the Kodaira--Spencer isomorphism $\Gr\nabla:\omega^{\otimes2}\xrightarrow{\sim}\Omega^1_X(\log D)$.
In particular, under the (logarithmic) $p$-adic Simpson correspondence defined by the canonical lift of $X$, $\mathbb{L}$ corresponds to $(E,\vartheta)$ \cite[Corollary~3.4.22]{DLLZ23}.

Choose compatible primitive roots of unity $\zeta_{Np^n}$, with
$\zeta_{Np^0}=\zeta_N$ and $\zeta_{Np^{n+1}}^p=\zeta_{Np^n}$.
Let $X_n=X_{Np^n}$ be the connected compactified modular curve of level $Np^n$
with Weil pairing $\zeta_{Np^n}$, and let $f_n:X_n\to X$ be the natural morphism. Then $X_n$ is defined over $K_n=K_0(\zeta_{p^n})$. Let $D_n$ be its reduced cusp divisor, the maps $(X_n,D_n)\to(X,D)$ are finite Kummer--\'etale. 
	The level structure trivializes the local system modulo $p^n$ on $(X_n,D_n)$:
\[
 f_n^*(\mathbb L^\circ/p^n\mathbb L^\circ)
 \simeq(\mathbf Z/p^n\mathbf Z)^{\oplus2}.
\]
Thus the pro-Kummer--\'etale tower $(X_n)_{n\ge 1}$ trivializes $\mathbb L$. 
\end{secnumber}

\begin{secnumber}\label{sss:modular-tower-obstruction}
	However, pulling back along the tower $(X_n)_{n\ge 1}$ does not make $\mathbb{L}$ into a small representation. Indeed, to apply the relative small $p$-adic Simpson correspondence \cite{AGrel} to $f_n^*\mathbb{L}$ for the canonical lifting of $X_n$ over $K_n$, the pullback $f_n^*\mathbb{L}$ must satisfy the smallness condition with respect to $K_n$.

Normalize $v_p$ by $v_p(p)=1$. Let $e_n=e(K_n/K_0)$,
let $\varpi_n$ be a uniformizer of $K_n$, and let
$\mathfrak d_{K_n/K_0}$ be the different of $K_n/K_0$. Set
\[
r_{K_n}=v_p(\varpi_n\mathfrak d_{K_n/K_0})=n-\frac{1}{p-1}+\frac{1}{e_n}.
\]
	Recall from \cite[Definition~3.4.19]{AGrel} that the relative smallness of a $\bC$-local system over $K_n$ requires an invariant $\oo$-lattice on which the monodromy action is trivial modulo $p^{\alpha}$ for some $\alpha>r_{K_n}+2/(p-1)$.  
	However, $f_n^*\mathbb{L}^{\circ}$ is only trivial modulo $p^n$, and is therefore \textit{not} small with respect to $K_n$, since
	\begin{equation}\label{eq:smallness-inequality-cover}
	 n<r_{K_n}+\frac{2}{p-1}
	   =n+\frac{1}{e_n}+\frac{1}{p-1}.
	\end{equation}

	This shows that the argument of reduction to the small equivalence discussed in \ref{sss:alternative-small} cannot be applied to $\mathbb{L}$ along the tower $(X_n)_{n\ge 1}$ with their canonical liftings simply by increasing $n$. A new idea is therefore needed to prove \Cref{t:intro-nilpotent} for the canonical lifting.
\end{secnumber}

\section{Hodge--Tate map on the pro-finite-\'etale cover trivializing the uniformizing representation}
\label{s:HT-map}
Let $C$ be a smooth projective curve over $\bC$ of genus $\ge 2$ and $\overline{x}$ a geometric point of $C$. 
Let $L$ be a line bundle over $C$ equipped with an isomorphism $\theta:L^{\otimes 2}\simeq \Omega_C$. 
By \Cref{t:intro-nilpotent}, there exists a lift $\mathcal{C}$ of $C$ over $\rB_{\dR,2}^+$ such that the uniformizing Higgs bundle $(E=L\oplus L^{-1},\Theta_{\theta})$, defined by $\theta$, corresponds to a continuous representation 
\begin{equation} \label{eq:rep-rho}
\rho: \pi_1(C,\overline{x})\to \GL_2(\bC)
\end{equation}
via the $p$-adic Simpson correspondence defined by $\mathcal{C}$. 

Choose a $\rho$-stable $\oo$-lattice and a basis, so that $\rho$ takes values in $\GL_2(\oo)$. For each $n\ge 1$, put
\[
 H_n=\Ker\bigl(\pi_1(C,\overline{x})
       \xrightarrow{\rho}\GL_2(\oo)
       \longrightarrow\GL_2(\oo/p^n\oo)\bigr),
\]
and let $C_n\to C$ be the connected pointed finite \'etale cover corresponding to $H_n$. 
In this section, we view $C$ and $C_n$ as adic spaces over $\Spa(\bC,\oo)$ via the analytification. 
Let $V_\rho$ be the pro-\'etale vector bundle associated to $\rho$ on $C$  \cite[\S~9.1]{HX24}. 
The pro-finite-\'etale cover of $C$ defined by the diamond
\begin{equation} \label{eq:cover-rho}
\widetilde{C}_{\rho}=\varprojlim_{n\ge 1} C_n
\end{equation}
trivializes the representation $\rho$ and $V_\rho$. 
Let $\pi:\widetilde C\to C$ be the natural map, which is equipped with a natural action of $\pi_1(C,\overline{x})$. 
On $\mathbb{P}^{1}_{\bC}$, viewed as the space of rank-one quotients
of $\bC^{\oplus2}$, has a natural action of $\GL_2$. 
\begin{theorem}
	\label{t:HT-map}	
	\textnormal{(i)}
	Assume in addition that $C$ descends to a smooth projective curve over $\overline K$. 
	The cover $\widetilde{C}_{\rho}$ is pointwise perfectoid in the sense of \cite{HeSen}, recalled in \ref{sss:pw-definition}. 

	\textnormal{(ii)} There exists a proper, surjective and equivariant Hodge--Tate map of diamonds over $\Spd(\bC,\oo)$:
	\[
	\pi_{\HT,\rho}:\widetilde{C}_{\rho}\to \mathbb{P}^{1,\dia}_{\bC}
	\]
	classifying the rank one quotient $\pi^*V_\rho\simeq \widehat{\mathscr{O}}^{\oplus 2}\twoheadrightarrow\pi^*\nu^*L$. 
	In particular, we have 
	\[\pi_{\HT,\rho}^*(\mathscr{O}_{\mathbb{P}^1}(1))\simeq \pi^*L.\] 
\end{theorem}

\begin{rem}
	In the case of modular/Shimura curve with canonical lift (see \S~\ref{ss:comments}), 
	the above construction is compatible with Scholze's definition of Hodge--Tate map \cite{Sch-HT}. 
\end{rem}

\begin{rem}\label{r:projective-uniformization} 
	Let $\overline{\rho}:\pi_1(C,\overline{x})\to \mathrm{PGL}_2(\bC)$ be the projective representation associated with $\rho$, and let $\widetilde{C}_{\overline{\rho}}\to C$ be the connected pointed pro-finite-\'etale cover corresponding to $\Ker(\overline{\rho})$.
	Since $\det(\rho)=1$, the quotient $\Ker(\overline{\rho})/\Ker(\rho)$ acts on $\bC^{\oplus 2}$ through $\{\pm I_2\}$ and hence trivially on $\mathbb{P}^{1,\dia}_{\bC}$.
	By equivariance, $\pi_{\HT,\rho}$ descends along $\widetilde{C}_{\rho}\to \widetilde{C}_{\overline{\rho}}$ to a proper, surjective, $\pi_1(C,\overline{x})$-equivariant map
	\[
	\pi_{\HT,\overline{\rho}}:\widetilde{C}_{\overline{\rho}}\longrightarrow \mathbb{P}^{1,\dia}_{\bC}.
	\]

	For a fixed lifting $\mathcal{C}$, one can check that both $\widetilde{C}_{\overline{\rho}}$ and $\pi_{\HT,\overline{\rho}}$ are independent of the choice of theta characteristic $L$ and depend only on $\mathcal{C}$.
\end{rem}
\subsection{Pointwise perfectoidness}\label{ss:pointwise-perfectoidness}

\begin{secnumber}\label{sss:pw-setup}
	Recall that the Higgs--Tate torsor $\mathcal{P}_{\mathcal{C}}$ under $\nu^*\xi T_C$ on $C_{\proet}$ associated with the lifting $\mathcal{C}$ \S~\ref{sss:Higgs-Tate-torsor} has Hodge--Tate image $\HT(\mathcal{P}_{\mathcal{C}})=1$ and defines an extension \Cref{p:Torsors-equiv}:
	\begin{equation}\label{eq:pw-lift-extension}
		0\to \xi \hatO_C\to \mathcal{F}_{\mathcal{C}}\to \nu^*(\Omega_C)\to 0.
	\end{equation}

	Recall that $V_\rho$ is the pro-\'etale vector bundle on $C_{\proet}$ associated with $\rho$ \eqref{eq:rep-rho}. 
	By \Cref{p:Torsors-equiv}, tensoring the short exact sequence $0\to \nu^*L^{-1}\to V_\rho\to \nu^*L\to 0$ with $\nu^*L$ and identifying $\theta:L^{\otimes 2}\simeq \Omega_C$ yields the extension \eqref{eq:pw-lift-extension}. 

	Since $V_\rho$ is trivial on $\widetilde{C}_{\rho}$, we obtain isomorphisms:
	\begin{equation}\label{eq:pw-twisted-triviality}
		\mathcal{P}_{\mathcal{C}}\simeq V_\rho\otimes \nu^* L,\qquad
	\pi^*\mathcal{P}_{\mathcal{C}}\simeq (\pi^*\nu^*L)^{\oplus 2}.
	\end{equation}
\end{secnumber}

\begin{secnumber}\label{sss:pw-descent}
	After enlarging $K$, we may assume there is a smooth projective curve $Y$ over $K$ such that $C=Y_{\bC}$. 
	Put $Z_0=Y\otimes_{K}\overline{K}$. The tower $(C_n)_{n\ge 1}$ over $C$ descends to a tower of connected finite \'etale covers $(Z_n)_{n\ge 1}$ of $Z_0$ over $\overline{K}$ such that $C_n\simeq Z_n\otimes_{\overline{K}}\bC$. 
	We set $Z_\infty=\varprojlim_{n\ge 1} Z_n$ and define the function fields
	\begin{equation}\label{eq:pw-function-field}
	\mathcal{K}=K(Y),\qquad \mathcal{L}=\varinjlim_{n\ge 0}\overline{K}(Z_n).
	\end{equation}
	Then $\mathcal{L}/\mathcal{K}$ is an algebraic extension containing $\overline{K}$, and $Z_\infty$ is the normalization of $Y$ in $\mathcal{L}$.
\end{secnumber}

\begin{secnumber}\label{sss:pw-definition}
	We recall the notion of \textit{pointwise perfectoidness} for the pro-finite \'etale cover $(Z_n)_{n\ge 1}$ in \cite{HeSen}. 
	A point $\widetilde{y}\in Z_\infty$ is a compatible system $(y_n)_{n\ge 1}$ with $y_n\in Z_n$, with residue field 
	\begin{equation}\label{eq:pw-residue-field}
		F:=\kappa(\widetilde{y})=\varinjlim_{n\ge 1}\kappa(y_n).
	\end{equation}

	Following \cite[Theorem 9.11(1)]{HeSen}, we say that $\widetilde{C}_{\rho}$ is \textit{pointwise perfectoid} if, for every point $\widetilde{y}\in Z_{\infty}$ and every height-one valuation ring $(W,w)\subset F$ extending the valuation on $\mathscr{O}_{\overline{K}}$ with $W[1/p]=F$, the valued field $(F,w)$ is pre-perfectoid. 
	Recall that this means that the valuation $w$ is non-discrete and the Frobenius $W/pW\to W/pW$ is surjective \cite[2.6]{HeSen}.

	Fix a point $\widetilde{y}\in Z_\infty$ and a height-one valuation ring $(W,w)\subset F=\kappa(\widetilde{y})$ as above.	
	Choose an algebraic closure $\overline{\mathcal{K}}$ of $\mathcal{L}$, a geometric point $\overline{y}=\Spec(\overline{F})$ in the normalization of $Y$ in $\overline{\mathcal{K}}$ lying over $\widetilde{y}\in Z_\infty$, and an extension of $w$ to $\overline{F}$. 
	Let $F^h\subset \overline{F}$ be the corresponding henselization and set
\[
H=\widehat{F^{h}},\qquad \Omega=\widehat{\overline{F}}.
\]
Since the valuation has height one, $F^h$ embeds into $\widehat{F}$ and $H=\widehat{F}$. 

Let $D\subset \Gal(\overline{\mathcal{K}}/\mathcal{L})$ be the valuative decomposition subgroup stabilizing $\overline{y}$ and its valuation (denoted by $\Gal(\overline{\mathcal{K}}/\mathcal{L}^{\mathrm{vh}})$ in \cite[\S~7.8]{HeSen}). 
	Its action on $\overline F$ has image $\Gal(\overline F/F^h)$.
	By the Ax--Sen--Tate theorem \cite[Lemma 7.18]{HeSen}, we have $\Omega^D=H$. 
\end{secnumber}

\begin{lemma}\label{l:pw-faltings-comparison}
	Let $\mathcal{P}_{\mathcal{C},\overline{y}}$ denote the $\Omega$-fiber of $\mathcal{P}_{\mathcal{C}}$ at $\overline{y}$. 
	Then there exists a $D$-equivariant isomorphism of $\Omega$-vector spaces
	\begin{equation}\label{eq:pw-faltings-comparison}
	\mathcal{P}_{\mathcal{C},\overline{y}}\simeq \mathscr{E}_{Y,\overline{y}},
	\end{equation}
	where $\mathscr{E}_{Y,\overline{y}}$ denotes the geometric Faltings extension of \cite[\S~6.4, \S~7.12]{HeSen}:
	\begin{equation}\label{eq:pw-he-extension}
		0\to \Omega(1)\xrightarrow{\iota} \mathscr{E}_{Y,\overline{y}}\to \Omega\otimes_{\mathscr{O}_{Y}}\Omega_{Y/K}\to 0.
	\end{equation}
\end{lemma}

\begin{proof} 
	We first consider the canonical lifting $\mathcal{C}^{\can}=(Y\otimes_{K}\rB_{\dR,2}^+)^{\ad}$.
	Choose an invertible \'etale coordinate $t$ on $Y$ around the image of $\overline{y}$ and a compatible system of $p^m$-th roots $(t^{1/p^m})_{m\ge 1}$ around $\overline{y}$. This allows us to define the cocycle
	\[
	c_t:D\to \Zp(1),\qquad \sigma\mapsto \bigl(\sigma(t^{1/p^m})/t^{1/p^m}\bigr)_{m\ge 1}.
	\]
	
	Let $d\log t$ be a $D$-invariant base of $\Omega\otimes_{\mathscr{O}_Y}\Omega_{Y/K}$. 
	Choose a lift $e$ of $d\log t$ in $\mathcal{P}_{\mathcal{C}}$ (or $\mathscr{E}_{Y,\overline{y}}$). 
	For any $\sigma\in D$, then $\sigma(e)-e\in \xi\Omega$. 
	The cocycle 
	\[
	D\to \xi \Omega,\qquad \sigma \mapsto \sigma(e)-e\in \xi\Omega
	\]
	determines the extension. It suffices to compare the cocycles defined by two extensions. 

	In the construction of $\mathscr{E}_{Y,\overline{y}}$ \cite[Theorems 7.13]{HeSen}, one choose a lift $e_t=(d\log(t^{1/p^m}))_{m\ge 1} \in \mathscr{E}_{Y,\overline{y}}$ of $d\log t$ has the cocycle $c_t(\sigma)$:
	\[
	\sigma(e_t)-e_t=(d\log \frac{\sigma(t^{1/p^m})}{t^{1/p^m}})=\iota (c_t(\sigma)).
	\]
	
	On the other hand, the local calculation of \cite[Lemma 2.16]{Heu25} identifies the cocycle of the Higgs--Tate torsor $\mathcal{P}_{\mathcal{C}^{\can}}$ with $c_t$. 
	This allows us to conclude the isomorphism \eqref{eq:pw-faltings-comparison}. 

	For a general lifting $\mathcal{C}$, the two liftings $\mathcal{C}$ and $\mathcal{C}^{\can}$ are isomorphic on every smooth affinoid $U\subset C$. Indeed, the torsor of isomorphisms lifting $\id_U$ is under $\xi T_U$, and $\rH^1(U,\xi T_U)=0$. 
	Choose such a neiborghood $U$ of the valuative point $\overline{y}$, and choose an isomorphism of liftings. 
	Functoriality identifies the extensions $\mathcal{P}_{\mathcal{C}}$ and $\mathcal{P}_{\mathcal{C}^{\can}}$ on $U_{\proet}$. 
	Combining it with the canonical case, we prove \eqref{eq:pw-faltings-comparison} for $\mathcal{C}$.
\end{proof}

\begin{proof}[Proof of \Cref{t:HT-map}(i)]
	Let $\widetilde{y}\in Z_\infty$ and a height-one valuation ring $W\subset\kappa(\widetilde{y})$ be as in \ref{sss:pw-definition}.
	Since $D\subset \Gal(\overline{\mathcal{K}}/\mathcal{L})$ fixes the chosen lift of $\overline{y}$ to each $Z_n$, its image in $\pi_1(C,\overline{x})$ lies in $\bigcap_{n\ge 1}H_n=\Ker(\rho)$.
	Consequently, the fiber $V_{\rho,\overline{y}}$ admits a $D$-invariant basis over $\Omega$.
	Moreover, the fiber of $\nu^*L$ at $\overline{y}$ comes by scalar extension from a one-dimensional $H$-vector space, and $D$ acts trivially on $\xi$.
	Combining \eqref{eq:pw-twisted-triviality} \Cref{l:pw-faltings-comparison}, we obtain an isomorphism of semilinear $D$-representations
	\begin{equation}\label{eq:pw-invariant-basis}
	\mathscr{E}_{Y,\overline{y}}\simeq \Omega\otimes_H H^{\oplus 2},
	\end{equation}
	where $D$ acts only through the first factor $\Omega$.
	Since $\Omega^D=\widehat{F^h}=H$, we have
	\begin{equation}\label{eq:pw-dimension}
	\dim_H\bigl(\mathscr{E}_{Y,\overline{y}}^D\bigr)=2.
	\end{equation}

	Since $\overline{K}\subset \mathcal{L}$ and $\dim Y=1$, it follows from \cite[Lemma 7.18, Theorem 7.19]{HeSen} that $(F,W)$ is pre-perfectoid. This finishes the proof.
\end{proof}

\subsection{Hodge--Tate map}

\begin{secnumber}
	Following Grothendieck \cite[Tags~01NE, 089V]{Stacks}, the projective space is represented by the functor:
	\[
	S\mapsto \{(\mathcal{M},q) | \mathcal{M} \textnormal{ line bundle on } S, q:\mathscr{O}_S^{\oplus 2}\twoheadrightarrow \mathcal{M}\}/\sim, 
	\]
	 where $\sim$ denotes isomorphism of quotients. 
	We need a similar description for the associated diamond $\mathbb{P}^{1,\dia}_{\bC}$.

	Let $T$ be a perfectoid space. 
	The categories of $\mathscr O_T$-vector bundles for the analytic and \'etale topologies and $\widehat{\mathscr O}_T$-vector bundles for the pro-\'etale topology are equivalent \cite[Theorem~3.5.8]{KL16}, \cite[Lemma~17.1.8]{SW20}. 
\end{secnumber}

\begin{prop}\label{l:P1-diamond}
	For every $T\in\Perf_{\bC}$, there is a bijection, functorial in $T$,
	\begin{equation}\label{eq:P1-functor}
	\mathbb{P}^{1,\dia}_{\bC}(T)\xrightarrow{\sim}\bigl\{(\mathcal{M},q)\bigm| \mathcal{M}\textnormal{ line bundle on }T,\ q:\mathscr{O}_T^{\oplus 2}\twoheadrightarrow\mathcal{M}\bigr\}/\sim,
	\end{equation}
	where $(\mathcal{M},q)\sim(\mathcal{M}',q')$ if there exists an isomorphism $\beta:\mathcal{M}\xrightarrow{\sim}\mathcal{M}'$ with $\beta\circ q=q'$.
	It sends a morphism $f:T\to\mathbb{P}^{1,\ad}_{\bC}$ to the pullback by $f$ of the universal quotient $\mathscr{O}^{\oplus 2}\twoheadrightarrow\mathscr{O}_{\mathbb{P}^1}(1)$.
\end{prop}

\begin{proof}
	We denote by $\Perf_{\bC}$ the category of perfectoid spaces over $\Spa(\bC,\oo)$. 
	By \cite[Definition~15.5]{Diamond}, we have $\mathbb{P}^{1,\dia}_{\bC}(T)=\Hom_{\Spa(\bC,\oo)}(T,\mathbb{P}^{1,\ad}_{\bC})$. 
	By the universal property of the analytification \cite[Proposition~3.8]{Hub94}, composition with the canonical morphism of locally ringed spaces $\mathbb{P}^{1,\ad}_{\bC}\to\mathbb{P}^1_{\bC}$ induces a bijection from this set to the set of morphisms of locally ringed spaces $(|T|,\mathscr{O}_T)\to\mathbb{P}^1_{\bC}$ over $\Spec(\bC)$.
	
	It suffices to show that, for any locally ringed space $(T,\mathscr{O}_T)$ over $\Spec(\bC)$, morphisms $T\to\mathbb{P}^1_{\bC}$ over $\Spec(\bC)$ correspond to pairs $(\mathcal{M},q)$ as in \eqref{eq:P1-functor}, where $\mathcal{M}$ is an invertible $\mathscr{O}_T$-module. 
	The proof follows the same argument in the case of schemes. 
	Given $f:T\to\mathbb{P}^1_{\bC}$, we take the pullback of the universal quotient. 
	Conversely, let $(\mathcal{M},q)$ be such a pair, and let $s_0,s_1\in\Gamma(T,\mathcal{M})$ be the images of the standard basis of $\mathscr{O}_T^{\oplus 2}$. 
	For $i=0,1$, let $T_i\subset T$ be the subset of points $x$ such that $s_i$ generates $\mathcal{M}_x$.
	Since $\mathcal{M}$ is invertible, $T_i$ is open and $s_i|_{T_i}$ is a basis of $\mathcal{M}|_{T_i}$.
	Since $q$ is surjective and $\mathscr{O}_{T,x}$ is local, Nakayama's lemma shows that $T=T_0\cup T_1$. 
	Let $t_0\in\Gamma(T_0,\mathscr{O}_T)$ (resp.\ $t_1\in\Gamma(T_1,\mathscr{O}_T)$) be the element such that $s_1=t_0s_0$ on $T_0$ (resp.\ $s_0=t_1s_1$ on $T_1$).
	By \cite[Tag~01I1]{Stacks}, $t_0$ defines a morphism $T_0\to D_+(x_0)=\Spec(\bC[x_1/x_0])$, and similarly for $t_1$.
	Since $t_0t_1=1$ on $T_0\cap T_1$, these morphisms glue to a morphism $T\to\mathbb{P}^1_{\bC}$.
	The two constructions are inverse to each other and compatible with isomorphisms of pairs. This proves the claim.
\end{proof}

\begin{secnumber}
	Recall that we have a short exact sequence of Higgs bundles on $C$:
	\[
	0\to (L^{-1},0)\to (E,\Theta)\to (L,0) \to 0.
	\]
	Applying the inverse $p$-adic Simpson correspondence $\rS_{\mathcal{C}}$ \eqref{eq:inverse-spectral}, defined by $\mathcal{C}$, we obtain a short exact sequence of pro-\'etale vector bundles on $C$:
	\[
	0\to \nu^*L^{-1} \to V_{\rho}\to \nu^*L\to 0,
	\]
	where $V_{\rho}$ is the pro-\'etale vector bundle associated with the representation $\rho$ \eqref{eq:rep-rho}. 
	
	The pullback of $V_{\rho}$ to the cover $\pi:\widetilde{C}_{\rho}\to C$ is trivial. 
	For any perfectoid space $T$ with a map $f:T\to \widetilde{C}_{\rho}$ of diamonds, we can define a quotient of rank two trivial pro-\'etale vector bundle:
	\begin{equation} \label{eq:pullback-quotient}
		\widehat{\mathscr{O}}_T^{\oplus 2}\simeq (\pi\circ f)^*_{\proet} V_{\rho}\twoheadrightarrow (\pi\circ f)^*_{\proet} (\nu^* L). 
\end{equation}
	Since $T$ is perfectoid, this quotient descends to a rank-one quotient of analytic $\mathscr{O}_T$-vector bundles. 
	This allows us to define a Hodge--Tate map relative to $\rho$:
	\[
	\pi_{\HT,\rho}:\widetilde{C}_{\rho} \to \mathbb{P}^{1,\dia}_{\bC},\qquad (f:T\to \widetilde{C}_{\rho})\mapsto \eqref{eq:pullback-quotient}.
	\]
\end{secnumber}

\begin{secnumber}
	We verify equivariance for the right actions. Set $V=\bC^{\oplus 2}$. 
For $f:T\to\widetilde C$, let
\[
 \tau_f:\widehat{\mathscr O}_T\otimes_{\bC}V
   \xrightarrow{\sim}(\pi\circ f)^*_{\proet}V_\rho
\]
be the trivialization defined by $f$. 
For every $g\in \pi_1(C,\overline{x})$, the action is given by
\[
 \tau_{f\cdot g}=\tau_f\circ(\id\otimes \rho(g)).
\]
Let $\varepsilon:V_\rho\twoheadrightarrow\nu^*L$ be the quotient
on $C$ and write
\[
 q_f=(\pi\circ f)^*\varepsilon\circ\tau_f: \widehat{\mathscr O}_T\otimes_{\bC}V \to (\pi \circ f)^*(\nu^*L). 
\]
for the quotient defining $\pi_{\HT,\rho}(f)$. Since
$\pi\circ(f\cdot g)=\pi\circ f$, the targets of the two quotients
are canonically the same, and we obtain
\[
 q_{f\cdot g}
 = (\pi\circ f)^*\varepsilon\circ\tau_{f\cdot g}
 = q_f\circ\rho(g).
\]
This shows that \[
 \pi_{\HT,\rho}(f\cdot g)
   =[q_f\circ\rho(g)]
   =\pi_{\HT,\rho}(f)\cdot g. 
\]
\end{secnumber}

\begin{secnumber}
	Next, we show the properness of $\pi_{\HT,\rho}$.  
	We write $\pi_{\HT}$ for $\pi_{\HT,\rho}$ for simplicity. 
	
	Let $S=\Spd(\bC,\oo)$. 
	The inverse limit $\widetilde{C}_{\rho}=\varprojlim C_n$ is qcqs over $S$ by \cite[Lemma 11.22]{Diamond}. 
	Each $C_n^{\dia}$ is proper over $S$. 
	In the valuative criterion, the unique extensions at the finite levels are compatible and therefore define a unique extension to $\widetilde C$. Thus \cite[Proposition 18.3]{Diamond} shows that $\widetilde{C}_{\rho}$ is proper over $S$. 
	The morphism $\pi_{\HT}$ decomposes as:
	\[
	\widetilde{C}_{\rho}\xrightarrow{\Gamma}\widetilde{C}_{\rho}\times_S \mathbb{P}^{1,\dia}_{\bC} \to \mathbb{P}^{1,\dia}_{\bC}.
	\]
	Since $\mathbb{P}^{1,\dia}_\bC$ is separated, we deduce that $\Gamma$ is a closed immersion. 
	And hence $\pi_{\HT}$ is proper. 
\end{secnumber}

\begin{secnumber}
	We first show the surjectivity of $\pi_{\HT}$ on classical points. 
	Consider the diagram:
	\[
	C \xleftarrow{\pi} \widetilde{C}_{\rho} \xrightarrow{\pi_{\HT}} \mathbb{P}^1_{\bC}.
	\]
	By construction, we have isomorphisms:
	\begin{equation} \label{eq:pullback-O(1)}
	\pi_{\HT}^*\mathscr{O}_{\mathbb{P}^1}(1)\simeq \pi^*(L),\qquad
	\pi_{\HT}^*\mathscr{O}_{\mathbb{P}^1}(2)\simeq \pi^*(\Omega_C).
\end{equation}

	Suppose a classical point $a\in \mathbb{P}^1(\bC)$ does not lie in the image of $\pi_{\HT}$. 
	Then $\mathscr{O}(2)$ is trivial on $(\mathbb{P}^1_{\bC}-a)^{\dia}\simeq \mathbb{A}^{1,\dia}_{\bC}$. 
	Therefore $\pi^*(\Omega_{C})$ is trivial on $\widetilde{C}_{\rho}$. 

	Let $q:C^{\univ}\to C$ be the universal pro-finite-\'etale cover of $C$, which dominates $\widetilde{C}_{\rho}$ \cite[Definition 9.1.1]{HX24}. 
	On the one hand, we have $\rH^0(C^{\univ}_{\proet},\widehat{\mathscr{O}})=\bC$ \cite[Proposition~3.9]{Heu24}. 
	On the other hand, we have $\dim \rH^0(C^{\univ}_{\proet},q_{\proet}^*\nu^*\Omega_C)=\infty$ by \cite[Lemma 9.3.8]{HX24}. 
	This contradicts to the fact $\pi^*(\Omega_C)$ is trivial on $\widetilde{C}_{\rho}$. Hence $\pi_{\HT}$ is surjective on the classical points. 

	Finally, since $\pi_{\HT}$ is proper, the image on the topological spaces $|\pi_{\HT}|:|\widetilde{C}_{\rho}|\to |\mathbb{P}^1_{\bC}|$ is closed. 
	The set of classical points $\mathbb{P}^1(\bC)$ is dense on $|\mathbb{P}^1_{\bC}|$. 
	Hence, $|\pi_{\HT}|$ is surjective. \hfill\qed 
\end{secnumber}

\end{document}